\documentclass{article}

 \usepackage[main, preprint]{neurips_2026}

\usepackage{microtype}
\usepackage{graphicx}
\usepackage{booktabs} %

\usepackage{hyperref}       %

\newcommand{\anonymousOrSubmission}[2]{#2}
\usepackage{placeins}

\usepackage[utf8]{inputenc} %
\usepackage[T1]{fontenc}    %
\usepackage{url}            %
\usepackage{amsfonts}       %
\usepackage{nicefrac}       %
\usepackage{xcolor}         %

\usepackage{amsmath,amsfonts}
\usepackage{array}
\usepackage[caption=false,font=normalsize,labelfont=sf,textfont=sf]{subfig}
\usepackage{textcomp}
\usepackage{stfloats}
\usepackage{url}
\usepackage{verbatim}
\usepackage[short]{optidef}
\usepackage{multicol}
\usepackage{multirow} 
\usepackage{comment}
\usepackage[shortlabels]{enumitem}

\usepackage{inconsolata}

\usepackage{threeparttable}
\usepackage[]{xcolor}
\usepackage{caption}
\usepackage{wrapfig}
\usepackage{booktabs}
\usepackage{pythonhighlight}
\usepackage{amsthm} %
\allowdisplaybreaks
\newtheoremstyle{newstyle}
{\parskip}%
{0pt}%
{}%
{0cm}%
{\bfseries}%
{.}%
{.5em}%
{\thmname{#1}\thmnumber{ #2}\thmnote{ (#3)}}
 
\theoremstyle{newstyle}

\newtheorem{assumption}{Assumption}%
\newtheorem{definition}{Definition}%
\newtheorem{theorem}{Theorem}%
\newtheorem{remark}{Remark}%
\newtheorem{example}{Example}%

\newcommand{\acados}{\texttt{acados}}

\newcommand{\grad}[1]{\nabla_{\!#1}}

\usepackage{mathrsfs, amsmath, amssymb, mathtools, amsfonts}
\DeclareMathSymbol{\shortminus}{\mathbin}{AMSa}{"39}
\usepackage{optidef}

\newcommand{\T}{^{\!\top}}

\newcommand{\inv}{^{\shortminus 1}}
\newcommand{\R}{\mathbb{R}}

\DeclareMathAlphabet{\mymathbb}{U}{bbold}{m}{n}

\newcommand{\ind}[1]{{_{\mathrm{#1}}}}
\newcommand{\nx}{n_x}
\newcommand{\nctrl}{n_u}

\newcommand{\Lag}{\mathcal{L}}
\DeclarePairedDelimiter\abs{\lvert}{\rvert}%
\DeclarePairedDelimiter\norm{\lVert}{\rVert}%

\makeatletter
\let\oldabs\abs
\def\abs{\@ifstar{\oldabs}{\oldabs*}}
\let\oldnorm\norm
\def\norm{\@ifstar{\oldnorm}{\oldnorm*}}
\makeatother

\newcommand{\eye}{%
	\text{\usefont{U}{bbold}{m}{n}1}%
}

\newcommand{\sol}{\mathrm{sol}}
\newcommand{\sadj}{s\ind{adj}}
\newcommand{\npd}{n_w}

\definecolor{Burgundy}{RGB}{144,0,32}

\title{Differentiable Nonlinear Model Predictive Control}

\author{%
	\textbf{Jonathan Frey}$^{1}$,
	\textbf{Katrin Baumgärtner}$^{1}$,
	\textbf{Gianluca Frison}$^{1}$,
	\textbf{Dirk Reinhardt}$^{2}$, \\
	\textbf{Jasper Hoffmann}$^{1}$,
	\textbf{Leonard Fichtner}$^{1}$,
	\textbf{Joschka Boedecker}$^{1}$,
	\textbf{Sebastien Gros}$^{2}$,
	\textbf{Moritz Diehl}$^{1}$ \\ \addlinespace[4pt]
	$^{1}$University of Freiburg \qquad
	$^{2}$NTNU Trondheim
}

\begin{document}
\maketitle

\begin{abstract}
The efficient computation of parametric solution sensitivities is a key challenge in the integration of learning-enhanced methods with nonlinear model predictive control, as their availability is crucial for many learning algorithms.
This paper discusses the computation of solution sensitivities of general nonlinear programs using the implicit function theorem and smoothed optimality conditions treated in interior-point methods (IPM).
We introduce sequential barrier quadratic programming (SBQP), which generalizes SQP with an IPM quadratic programming (QP) solver by allowing the barrier parameter to approach a prescribed nonnegative target value, recovering classical SQP when this target is zero.
Previous works are limited to convex or unconstrained formulations, or lack an implementation for efficient sensitivity evaluation.
The publication is accompanied by an efficient open-source implementation within the \acados{} framework, providing both forward and adjoint sensitivities for general optimal control problems, achieving speedups exceeding 3x over the state-of-the-art solvers \texttt{mpc.pytorch} and \texttt{cvxpygen}.
\end{abstract}

\section{Introduction}

In recent years, great research efforts have been made to combine the paradigms of learning and optimal control -- in particular nonlinear model predictive control (MPC) -- exploiting their complementary advantages.
The outcomes of these efforts cover various fields: from MPC with learned system models to imitation learning and behavior cloning to MPC-based Reinforcement Learning (RL) \citep{Goerges2017, Hewing2020a, Reiter2026}.
In particular, the approach of performing learning-enhanced adaptations of MPC schemes is very general and promises to deliver excellent performance while attaining desirable properties of MPC schemes, such as (robust) constraint satisfaction and stability \citep{Gros2020, Gros2022, Kordabad2022,Kordabad2022a, Kordabad2024}.
The approach has been successfully employed for a wide range of applications such as
car racing~\citep{Reiter2025a, Jahncke2026},
quadcopter control~\citep{Romero2024},
energy-management systems in residential buildings~\citep{Cai2023},
microgrid operation for demand response~\citep{Cai2023a}, and greenhouse climate control~\citep{Mallick2025}.
A wider adoption of this approach in control practice has been limited by the availability of dedicated software and, in particular, efficient solvers that provide parametric solution sensitivities -- a crucial prerequisite for a variety of learning algorithms \citep{Reinhardt2024a}.
Inherent constraint satisfaction is a key advantage of learning-enhanced MPC over purely data-driven approaches.
However, existing differentiable MPC implementations lack support for inequality constraints (cf. Table~\ref{tab:diff_mpc}), which has so far prevented real-world applications from benefiting from this advantage \citep{Rickenbach2025}.
Other techniques that can profit from OCP solution sensitivities include Sobolev learning for approximate MPC \citep{Lueken2023} and differentiable robust MPC \citep{Oshin2024}. %

As nonlinear MPC heavily relies on the numerical solution of nonlinear optimization problems in real time, tailored algorithms and efficient implementations have been a major research focus \citep{Frison2013a, Vanroye2023, Wang2010a, Zavala2009, Ferreau2014, Frey2020, Stellato2020, Arnstrom2022, Goulart2024, Ferreau2017, Houska2011e, Zanelli2017b}. %
In particular, the open-source software framework \texttt{acados}~\citep{Verschueren2021} implements many of these algorithms.
The numerical solvers provided by \texttt{acados} are tailored to the problem structure specific to optimal control, thus achieving exceptional computational efficiency, while keeping the problem formulation as general as possible, covering general nonlinear and parametric costs, dynamics and inequality constraints, which can be different at each stage of the problem~\citep{Frey2024b}.

The focus of this paper is the efficient computation of parametric solution sensitivities in the context of optimal control, as their availability -- and efficient computation -- is a key building block for a successful integration of learning-enhanced methods and optimization-based control.

\textbf{Contribution.}
The contribution of this paper is threefold:
First, this work discusses the efficient computation of parametric (adjoint) solution sensitivities for general nonlinear programs (NLP) via the implicit function theorem (IFT).
In particular, this work highlights common pitfalls of the IFT and discusses their mitigation.
Previous work in the machine learning (ML) literature was subject to this pitfall, or limited to convex or unconstrained nonlinear least squares problems.
Second, we introduce SBQP, which generalizes SQP with an IPM QP solver by allowing the barrier parameter to approach a prescribed nonnegative target value.
We detail how IPMs or SBQP can be used to obtain a differentiable approximation of the solution map and its sensitivities suitable for gradient-based learning algorithms by leveraging the smoothing of the optimality conditions that these methods employ.
Finally, we present an efficient implementation of a differentiable solver for optimal control structured NLPs within the open-source software package \acados{} which provides both forward and adjoint sensitivities.
In particular, the optimal control problem (OCP) formulation covers general, potentially nonlinear, costs and constraints on both states and controls while allowing all problem functions, i.e., costs, constraints and dynamics, to be fully parametric.
A comparison with state-of-the-art solvers demonstrates speedup factors greater than three on a CPU.

\textbf{Related work.}
One approach to analyze and compute parametric solution sensitivities of general nonlinear programs is to apply the IFT~\citep{Dini1907} to the necessary conditions of optimality, also called Karush-Kuhn-Tucker (KKT) conditions~\citep{Karush1939, Nocedal2006}.
An alternative approach differentiates through the iterations of an optimization algorithm~\citep{Tamar2017}.
In the following, we adopt the former approach, which offers superior memory efficiency and computational complexity compared to differentiating through the iterations of an optimization solver, as demonstrated in the context of MPC by \citet{Amos2018} using iLQR.

A qualitative comparison of available software tools for differentiable MPC is given in Table~\ref{tab:diff_mpc} and important relations are summarized next.
\citet{Pirnay2012} details solution sensitivity computation via the smoothed KKT system tackled in an IPM, providing useful theoretical results and the implementation \texttt{sIPOPT}.
This implementation builds on the popular NLP solver \texttt{IPOPT}~\citep{Waechter2006}, but requires defining parameters as variables and does not offer adjoint sensitivities.
In~\citet{Andersson2018a}, the computation of forward and adjoint solution sensitivities based on an active-set QP solver was presented.
More recently, the analysis of parametric solution sensitivities has received a lot of attention in the ML community.
\citet{Amos2017} suggested embedding the solution of an optimization problem in a neural network, but their derivations and implementation are limited to quadratic programs (QPs).
In \citet{Lee2019}, the differentiable QP solver from \citet{Amos2017} was leveraged in the context of meta-learning.
In \citet{Agrawal2019}, neural network architectures comprising so-called differentiable convex optimization layers are introduced.
These layers contain an optimization problem formulated via disciplined parametrized programming, closely related to disciplined convex programming~\citep{Grant2006} and are implemented in~\texttt{cvxpylayers}.
Recently, \texttt{cvxpygen}~\citep{Schaller2022} was extended to support solution sensitivities which showed significant speedups compared to the implementation in \texttt{cvxpylayers}~\citep{Schaller2026}.
The software package \texttt{Theseus} provides solution sensitivities for nonlinear least-squares optimization problems~\citep{Pineda2022}, but is limited to soft constraints.
\citet{Mandi2020} solve linear programs with an interior-point approach and stop the barrier parameter at a certain threshold to differentiate the solution of a smoothed KKT system.
However, their work does not consider an efficient backward pass via adjoint solution sensitivities.
In contrast to implicit smoothing via the IPM, the nondifferentiability of the solution map might be directly tackled using a nonsmooth variant of the IFT which is based on the assumption of path differentiability and the concept of conservative Jacobians \citep{Davis2020, Bolte2021}.

While the works mentioned above consider general, i.e. structure-agnostic, solvers, \citet{Amos2018} first considered differentiable solution maps within the context of structure-exploiting solvers tailored to OCPs coining the term \textit{differentiable MPC}.
However, their approach is restricted to problems with quadratic cost functions and simple bounds on the controls and only covers parametric nonlinear dynamics and parametric quadratic costs. %
\citet{Adabag2026} present an implementation of differentiable MPC for GPU, which tackles a more restrictive formulation with no inequalities and a cost function that is decoupled in controls and states.
The papers by \citet{Amos2018, Adabag2026} also suggest the computation of solution sensitivities via the KKT system of a convex approximation of the problem, an approach which might yield highly degraded sensitivity results as found in \citet{Dinev2022} and shown in the present paper, see Remark~\ref{rem:ift_pitfall} and Example~\ref{example:solution_sens_pendulum}.
\citet{Jin2020} suggest computing solution sensitivities via Pontryagin's Maximum Principle, but the approach limited to unconstrained optimal control.
\citet{Bounou2023} suggest a proximal solution approach to obtain sensitivities for OCPs.
All works cited in this paragraph lack the ability to handle inequality constraints which are parametric, nonlinear or include the states.
In contrast, \citet{Zuliani2025} use the nonsmooth IFT variant to compute sensitivities in an MPC context considering convex quadratic programs, but do not provide an efficient implementation exploiting the OCP structure.

\section{Preliminaries} %
\label{sec:preliminaries}
This section briefly reviews the mathematical foundations and formal requirements underlying the sensitivity analysis of general parametric nonlinear programs (NLP) via the implicit function theorem (IFT), stated in~\ref{sec:ift}.
Before diving into the mathematical details, we introduce a tutorial example which illustrates nondifferentiability of the solution map -- a property that might pose major challenges in practice, in particular in the context of gradient-based learning.

\begin{figure}
	\centering
	\includegraphics[width=\columnwidth]{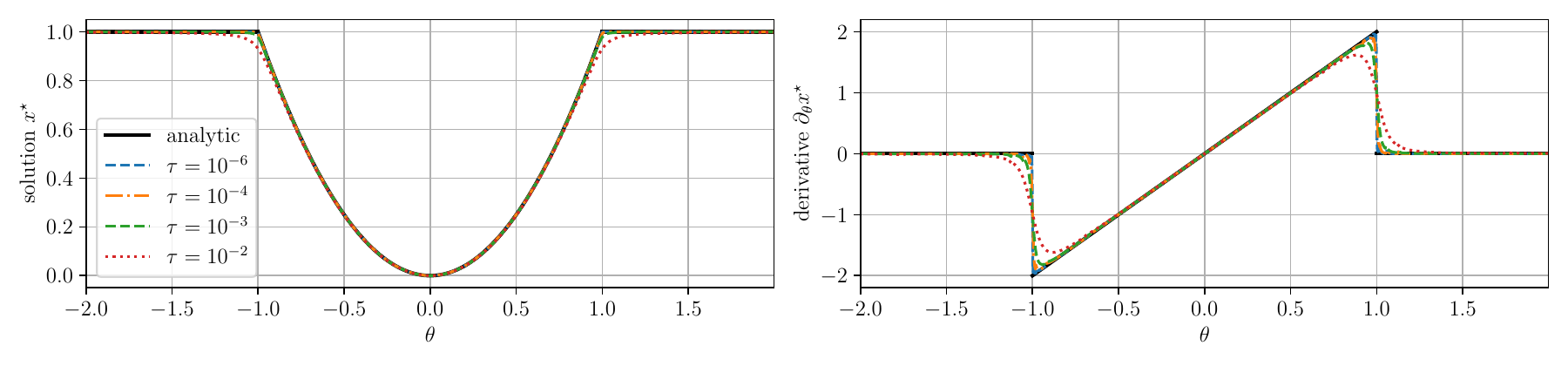}
	\vspace{-6mm}
	\caption{Simple example visualizing a nonsmooth solution map, its derivative and some smooth approximations.
	\label{fig:solution_sens_non_ocp}}
	\vspace{-5mm}
\end{figure}

\begin{example}\label{example:solution_sens_tutorial}
We regard the following scalar parametric optimization problem:
\begin{equation}
\min_x ~(x-\theta^2)^2\; \text{ s.t. } \; -1 \leq x \leq 1. \label{eq:simple_parametric_nlp}
\end{equation}
\vspace{-1.5mm}

The solution map and its sensitivity are given by
\vspace{-2mm}
\begin{align*}
x^\star(\theta) &= \begin{cases}
\theta^2, \; \textrm{if} \; \theta \in [-1, 1] \\
1, \; \textrm{otherwise},
\end{cases}
\quad
\partial_\theta x^\star(\theta) = \begin{cases}
2\cdot \theta, \; &\textrm{if} \; \theta \in (-1, 1) \\
0, \; &\textrm{if} \; \abs{\theta} > 1 \\
\textrm{not defined}, \; &\textrm{if} \; \theta \in \{-1, 1\}.
\end{cases}
\end{align*}
\end{example}
Figure~\ref{fig:solution_sens_non_ocp} visualizes the solution map as well as its derivative and their smooth approximations, which will be explained later in Section \ref{sec:dense_solution_sens}.

\subsection{Solving nonlinear programs with an interior-point based SQP method}
\label{sec:sqp-ipm}
We are interested in differentiating the solution of a nonlinear program of the form
\begin{align}\label{eq:z_star_def}
z^\sol(\theta) \coloneqq \arg \min_{z \in \R^{n_z}} f(z; \theta)\; \text{ s. t. } \;\; g(z; \theta) = 0,~ h(z; \theta) \leq 0.
\end{align}
We assume that the objective function $f:\R^{n_z} \times \R^{n_\theta} \rightarrow \R$ and the constraint functions $ g: \R^{n_z}\times\R^{n_\theta} \to \R^{n_g}$, $ h: \R^{n_z}\times\R^{n_\theta} \to \R^{n_h}$ are twice continuously differentiable jointly in $(z,\theta)$.
Introducing the Lagrange multipliers $\lambda \in \R^{n_g}$ and $\mu \in \R^{n_h}$, the Lagrangian function for \eqref{eq:z_star_def} is given by
	\vspace{-2mm}
\begin{align}
	\Lag(z, \lambda, \mu; \theta) = f(z; \theta) + \lambda\T g(z; \theta) + \mu \T h(z; \theta).
\end{align}
The Karush-Kuhn-Tucker (KKT) conditions of NLP~\eqref{eq:z_star_def} can be written as

\vspace{-2mm}
\begin{subequations}
	\noindent \centering
	\label{eq:nlp_kkt}
	\begin{minipage}{0.6\textwidth}
	\begin{align}
	\grad{z}  f(z; \theta) + \grad{z} g(z; \theta) \lambda &+ \grad{z} h(z; \theta) \mu = 0, \label{eq:KKT_stat} \\
	g(z; \theta) &= 0, \\
	h(z; \theta) &\leq 0,
	\end{align}
	\end{minipage}
	\hfill
	\begin{minipage}{0.399\textwidth}
		\begin{align}
			\mu & \geq 0, \label{eq:KKT_mult_pos} \\
			\mu_i h_i(z; \theta) &= 0, \, i=1,\dots,n_h. \label{eq:KKT_comp}
		\end{align}
	\end{minipage}\bigskip
\end{subequations}

If there exists $\lambda^\star, \mu^\star$ such that $v^\star \!=\! (z^\star, \lambda^\star, \mu^\star)$ is a solution of \eqref{eq:nlp_kkt}, $v^\star$ is called a KKT point.
A KKT point is a local minimizer of~\eqref{eq:z_star_def} if the following assumption holds, cf.~\citet{Nocedal2006}.
\begin{assumption} \label{as:regularity}
	Let $v^\star$ be a KKT point of the NLP \eqref{eq:z_star_def} satisfying linear independence constraint qualification (LICQ), second-order sufficient conditions of optimality (SOSC) and strict complementarity.
\end{assumption}

\textbf{Sequential quadratic programming.}
Sequential quadratic programming (SQP) is an efficient solution method for solving NLPs of the form \eqref{eq:z_star_def}.
In the context of nonlinear MPC, SQP is particularly popular due to its warm-starting capabilities and real-time variants~\citep{Diehl2005, Nurkanovic2019a, Frey2024a}.
For a given primal-dual initialization $(z_0, \lambda_0, \mu_0)$ and using the compact notation $v_k = (z_k, \lambda_k, \mu_k)$, an SQP method proceeds by successively solving quadratic programs (QPs) of the form
\begin{align}
	\begin{split}
\label{eq:dense_QP}
\Delta z_\mathrm{QP}^\sol(\theta, v_k) \coloneqq ~\arg \min_{\Delta z} ~q_k\T \Delta z + \tfrac{1}{2} \Delta z \T Q_k \Delta z
\quad \text{ s.t. } \quad
g_k + G_k \Delta z = 0, \
h_k + H_k \Delta z \leq 0,
	\end{split}
\end{align}
where $Q_k$ is an approximation of the Hessian of the Lagrangian $\nabla^2_{zz}\mathcal{L}(z_k, \lambda_k, \mu_k; \theta)$, the vector $q_k = \grad{z} f(z_k; \theta)$ is the objective gradient, $g_k = g(z_k; \theta)$ and $h_k = h(z_k; \theta)$ are the constraint evaluations and $G_k = \nabla g(z_k; \theta)^\top $, $H_k = \nabla h(z_k; \theta)^\top $ the Jacobians of the constraints.
The full-step SQP method updates its iterates by
\begin{align}
z_{k+1} = z_k + \Delta z_\mathrm{QP}^\sol(\theta, v_k), \quad
\lambda_{k+1} = \lambda_{\mathrm{QP}}^\sol(\theta, v_k), \quad
\mu_{k+1} &= \mu_{\mathrm{QP}}^\sol(\theta, v_k),
\end{align}
where %
$\mu_{\mathrm{QP}}^\sol(\theta, v_k)$, $\lambda_{\mathrm{QP}}^\sol(\theta, v_k)$ denote the dual solution of QP~\eqref{eq:dense_QP}.
The QP subproblems can be tackled with different solution strategies, which can be classified into active-set methods, interior-point methods and first-order methods.
This paper focuses on interior-point methods for solving the QPs, which are detailed in the following.

\begin{remark}[Hessian approximations]%
	\label{remark:hess_approx}
Many SQP approaches use an approximate Hessian $Q_k$ instead of the exact Hessian of the Lagrangian.
A popular class of Hessian approximations are the Gauss-Newton Hessian approximation and its generalizations.
Their advantages are that they are inherently positive-semidefinite, and represent outer convex structure in the subproblems which leads to beneficial convergence properties \citep{Bock2007b, Messerer2021a, Baumgaertner2022, Baumgaertner2023a}.
Additionally, these Hessian approximations are cheaper to compute compared to the exact Hessian.
In the context of unconstrained optimal control, the Gauss-Newton Hessian is used, for instance by the popular iLQR method \citep{Li2004a, Todorov2005, Amos2018}. %
\end{remark}

\textbf{QP solution via interior-point methods.}
Within an interior-point QP solver, the following system of smoothed KKT conditions of QP \eqref{eq:dense_QP} is solved:

\vspace{-3mm}
\begin{subequations}
\label{eq:ipm_qp_kkt}
\begin{minipage}[t]{0.56\textwidth}
\vspace{0pt}
\vspace{-\abovedisplayskip}
\begin{align}
    Q_k \Delta z + q_k + G_k\T \lambda + H_k\T \mu &= 0, \label{eq:ipm_stat} \\
    g_k + G_k \Delta z &= 0, \label{eq:ipm_eq} \\
    h_k + H_k \Delta z + s &= 0, \label{eq:ipm_ineq}
\end{align}
\end{minipage}
\hfill
\begin{minipage}[t]{0.42\textwidth}
\vspace{0pt}
\vspace{-\abovedisplayskip}
\begin{align}
    \mu_i s_i &= \tau, \;\, i=1,\dots, n_h, \label{eq:ipm_comp} \\
    s, \mu &> 0, \label{eq:ipm_positivity}
\end{align}
\end{minipage}
\end{subequations}
\vspace{2mm}

where the slack vector $s\in \R^{n_h}$ and the barrier parameter $\tau>0$ are introduced.
This system of equations is solved repeatedly for values of $\tau$ that approach zero from above.
In an IPM, the positivity conditions~\eqref{eq:ipm_positivity} are enforced with strict inequality, typically by applying some fraction to the boundary rule.
An IPM applies Newton steps to solve the system of equations consisting of~\eqref{eq:ipm_stat} - \eqref{eq:ipm_comp}.
The linear system solved in an interior-point iteration $j$ with barrier parameter $\tau_j$ can be written as
\vspace{-2mm}
\begin{align}
	\label{eq:ipm_ls}
	\underbrace{
	\begin{bmatrix}
		Q_k & G_k\T & H_k\T  & 0 \\
		G_k & 0 & 0 & 0 \\
		H_k & 0 & 0 & \eye \\
		0 & 0 & S_j & M_j
	\end{bmatrix}
	}_
	{\eqqcolon \mathcal{M}_k(s_j, \mu_j; \theta)}
	\begin{bmatrix}
		\Delta \hat{z} \\
		\Delta \hat{\lambda} \\
		\Delta \hat{\mu} \\
		\Delta \hat{s}
	\end{bmatrix}
	= \!\!\!\!\!\!\!\!\!\!\!
	\underbrace{-
	\begin{bmatrix}
		\hat{q}_j \\
		\hat{g}_j \\
		\hat{h}_j \\
		\hat{m}_j
	\end{bmatrix}
	}_{\eqqcolon r_k(\Delta z_j, \lambda_j, \mu_j, s_j; \tau_j, \theta)}
\end{align}
with diagonal matrices $S_j = \mathrm{diag}(s_j)$ and $ M_j = \mathrm{diag}(\mu_j) $ and residuals
\begin{align*}
	\hat{q}_j &= Q_k \Delta z_j \!+\! q_k \!+\! G_k\T \lambda_j \!+\! H_k\T \mu_j, \;\;
	\hat{g}_j = G_k \Delta z_j \!+\! g_k, \\
	\hat{h}_j &= H_k \Delta z_j \!+\! h_k \!+\! s_j, \;\;
	\hat{m}_j = (\mu_{j,i} s_{j,i} - \tau_j)_{i=1,\dots, n_h}, %
\end{align*}
where $(\Delta z_j, \lambda_j, \mu_j, s_j)$ denotes the current iterate of the IPM.
The linear system \eqref{eq:ipm_ls} can be reduced by eliminating $\Delta \hat{s}$ and $\Delta \hat{\mu}$ to arrive at a lower dimensional system with symmetric coefficient matrix
\begin{align}
	\label{eq:m_tilde}
	\widetilde{\mathcal{M}}_k(s_j, \mu_j; \theta) =
		\begin{bmatrix}
			Q_k + H_k\T S_j\inv M_j H_k & G_k\T \\
			G_k & 0
		\end{bmatrix}.
\end{align}

Eliminating $\Delta \hat{s}$ and $\Delta \hat{\mu}$ in this way is especially beneficial in the context of OCPs, as the reduced system can be solved very efficiently with a Riccati factorization, see Sec.~\ref{sec:diff_mpc}.
Invertibility of $\mathcal{M}_k$ follows from Assumption~\ref{as:regularity} for sufficiently small~$\tau$, see Thm.~\ref{thm:smooth_KKT_sol},~\ref{thm:smooth_KKT_item2}.

\textbf{Sequential barrier quadratic programming (SBQP).}
Next, we detail SBQP, which generalizes SQP with an interior-point QP solver by allowing the barrier parameter to approach a prescribed target $\tau\ind{min} \geq 0$ rather than zero, such that classical SQP is recovered when $\tau\ind{min} = 0$.
When $\tau\ind{min} > 0$, each SBQP iteration corresponds to solving the KKT conditions of QP \eqref{eq:dense_QP}, with the inequalities replaced by logarithmic barrier terms
\vspace{-.2cm}
\begin{align}
	\begin{split}
		\label{eq:bqp}
		\min_{\Delta z, s} ~q_k\T \Delta z + \tfrac{1}{2} \Delta z \T Q_k \Delta z  \shortminus \tau\ind{min} \sum_{i=1}^{n_h} \log(s_i)
		\ \text{ s.t. } \
		g_k + G_k \Delta z & = 0, \\[-5pt]
		h_k + H_k \Delta z & + s = 0.
	\end{split}
\end{align}
We refer to \eqref{eq:bqp} as a \textit{barrier QP} and denote its solution map by $(\Delta z_\mathrm{BQP}^\sol(\tau\ind{min}; \theta, v_k), s_\mathrm{BQP}^\sol(\tau\ind{min}; \theta, v_k))$.
In an IPM QP solver, the barrier parameter is iteratively reduced until a prescribed tolerance on the KKT residuals is met.
In SBQP, instead, we let $\tau_j$ converge to the target value $\tau\ind{min}$ from above and check the IPM residuals \eqref{eq:ipm_qp_kkt} with $\tau=\tau\ind{min}$.
For $\tau\ind{min} > 0$, SBQP corresponds to solving the nonlinear smoothed interior-point KKT system~\eqref{eq:nlp_kkt_ipm}, which consists of \eqref{eq:KKT_stat}-\eqref{eq:KKT_mult_pos} and instead of the complementarity condition \eqref{eq:KKT_comp}, it contains
\begin{align}
\mu_i h_i(z; \theta) &= -\tau\ind{min}, \, i=1,\dots,n_h. \label{eq:KKT_comp_ipm}
\end{align}
Note that the termination criterion of the SBQP solver needs to be adjusted accordingly with respect to SQP.
This allows us to obtain smoothed approximations of the solution map with valuable properties, see Thm.~\ref{thm:smooth_KKT_sol}.
These approximations are visualized for different $\tau\ind{min}$ values in Fig.~\ref{fig:solution_sens_non_ocp} and Fig.~\ref{fig:solution_sens_pendulum}.
Appendix~\ref{sec:ipm_sqp_details} provides more details on the algorithm presented in this section.

\subsection{Parametric solution sensitivities within SBQP}
\label{sec:dense_solution_sens}
This section discusses the computation of (approximate) parametric solution sensitivities of NLP~\eqref{eq:z_star_def} within SBQP.
On the theoretical side, two levels need to be considered:
First, we regard the QP~\eqref{eq:dense_QP} and review the conditions under which its solution sensitivities match those of the NLP.
Second, we review the conditions under which the smoothed KKT system~\eqref{eq:ipm_ls} yields an approximate solution to the QP and discuss its parametric sensitivities.

\begin{theorem}[NLP solution sensitivity existence] \label{thm:nlp-sens-exist}
Suppose Assumption~\ref{as:regularity} holds at a KKT point $z^\star$ of \eqref{eq:z_star_def} with parameter $\bar \theta$.
In a neighborhood of $\bar \theta$, there exists a differentiable function $z^\sol(\theta)$ with $z^\sol(\bar \theta) = z^\star$  that corresponds to a locally unique solution of \eqref{eq:z_star_def}.
\end{theorem}

Theorem~\ref{thm:nlp-sens-exist} is proved via the IFT in~\ref{sec:solution_sens_nlp}.
Next, we consider the smoothed KKT system \eqref{eq:nlp_kkt_ipm} as solved within an IPM or SBQP.
The following theorem shows how the solution map of NLP \eqref{eq:z_star_def} can be obtained as the solution to the smoothed KKT system as $\tau \rightarrow 0$.

\begin{theorem}[Smoothed KKT system]
	\label{thm:smooth_KKT_sol}
Suppose Ass.~\ref{as:regularity} holds at a KKT point $z^\star$ of NLP \eqref{eq:z_star_def} with multipliers $\lambda^\star$, $\mu^\star$.
Then, for all sufficiently small $\tau \coloneqq \tau\ind{min} > 0$, the following holds:
\begin{enumerate}[(i), noitemsep, nosep, leftmargin=20pt]
	\item There exists a solution of the smoothed interior-point KKT system $z^\sol_\textsc{ipm}(\tau; \bar{\theta})$ which is a continuously differentiable function with $ \lim_{\tau\rightarrow 0^+} z^\sol_\textsc{ipm}(\tau; \bar{\theta}) = z^\sol(\bar{\theta}) $ and $\norm{z^\sol_\textsc{ipm}(\tau; \bar{\theta}) - z^\star} \in \mathcal{O}(\tau)$.
	\item In a neighborhood of $\bar{\theta}$, there exists a differentiable function
$v(\tau; \theta) = (z(\tau; \theta), \lambda(\tau; \theta), \mu(\tau; \theta))$
	that corresponds to a locally unique solution of the smoothed interior-point KKT system \eqref{eq:nlp_kkt_ipm} and $v(0; \bar{\theta}) \coloneqq \lim_{\tau\rightarrow 0^+} v(\tau; \bar{\theta}) = (z^\star, \lambda^\star, \mu^\star)$ holds. \label{thm:smooth_KKT_item2}
\end{enumerate}
\end{theorem}
Theorem~\ref{thm:smooth_KKT_sol} can be deduced from Prop.~1 and~3 in \citet{Pirnay2012}.
A comprehensive proof is accompanying Theorem 8 in \citet{Wright1992}.
The next theorem states that the solution sensitivities of the smoothed interior-point KKT systems match the ones of the SBQP subproblem at convergence.

\begin{theorem}[Barrier NLP and barrier QP sensitivities coincide]
\label{thm:barrier-bqp}
Fix $\tau>0$ and a parameter $\bar{\theta}$, and let $w^\star = (z^\star, s^\star, \lambda^\star, \mu^\star)$ be a KKT point of the barrier NLP~\eqref{eq:dense_nlp_log_barrier} with $v^\star \coloneqq (z^\star, \lambda^\star, \mu^\star)$; assume LICQ and SOSC hold at $w^\star$.
Regard the barrier QP~\eqref{eq:bqp} obtained at the KKT point, and assume the exact Hessian $Q_\star = \nabla^2_{zz}\mathcal{L}(z^\star, \lambda^\star, \mu^\star; \bar\theta)$ is used.
Denote the solution map of the barrier NLP as $(z^\sol_\textsc{ipm}, s^\sol_\textsc{ipm})$ and the solution map of the barrier QP as $(\Delta z^\sol_\mathrm{BQP}(\tau; \theta, v^\star), s^\sol_\mathrm{BQP}(\tau; \theta, v^\star))$.
Then, the solution maps $z^\sol_\textsc{ipm}(\tau; \theta)$ and $z^\star + \Delta z^\sol_\mathrm{BQP}(\tau; \theta, v^\star)$, and their sensitivities, $\frac{\partial z^\sol_\textsc{ipm}}{\partial \theta}(\tau; \theta)$ and $\frac{\partial \Delta z_\mathrm{BQP}^\sol}{\partial \theta}(\tau; \theta, v^\star)$ coincide at $\bar{\theta}$.
The same holds for the solution maps of the Lagrange multipliers.
\end{theorem}
The proof is given in Appendix~\ref{sec:bqp_sens_match}.
Similarly, it is shown in Appendix~\ref{sec:qp_sens_match} that the QP sensitivities match the ones of the NLP, which corresponds to the proposed method for $\tau\ind{min}=0$, i.e. when SBQP corresponds to SQP.

Solutions to the smoothed KKT systems can be of particular interest in a learning context.
First, their smoothness can be beneficial, as shown in an RL context in \citet{Kordabad2021}.
Second, their computation can be computationally cheaper, as less interior point iterations might be needed to shrink the barrier parameter.
\begin{remark}[Differentiating across active set changes]
If Ass.~\ref{as:regularity} is not satisfied, the solution map might be nondifferentiable, as in Example~\ref{example:solution_sens_tutorial}, or even exhibit jumps, as shown by the example in~\ref{sec:jump_ocp}.
In case that only strict complementarity is violated in Ass.~\ref{as:regularity}, but strong SOSC holds, the proposed method provides the correct sensitivities of the smoothed solution map, as shown in~\ref{app:reg_log_barrier}.
\end{remark}

\textbf{Forward and adjoint sensitivity computation.}\label{sec:adjoint_sens}
In the remainder of this section, we use the compact notation $w \coloneqq (z, \lambda, \mu, s)\in \R^{\npd}$ with $n_w = n_z + n_g + 2n_h$ and assume that SBQP has converged to a solution $w^\star$ of the barrier KKT system \eqref{eq:nlp_kkt_ipm} for a fixed $\tau>0$.
We assume the conditions of Theorem~\ref{thm:barrier-bqp}, under which the barrier QP sensitivities coincide at $\bar{\theta}$ with those of the smoothed NLP.
We compute them from the barrier QP KKT system but denote them by the sensitivities of $w^\sol_\textsc{ipm}$, the quantity of interest.
Applying the IFT to the barrier QP KKT system gives
\vspace{-0.12cm}
\begin{align}
	\label{eq:forward_ipm_sens}
	\tfrac{\partial w^\sol_\textsc{ipm}}{\partial \theta}(w^\star; \tau, \theta) = -\mathcal{M}_\star(w^\star; \tau, \theta)\inv J_\star(w^\star; \tau, \theta),
\end{align}
where we introduced $J_\star(\cdot) \coloneqq \frac{\partial r_\star}{\partial \theta}(\cdot)$ via the function $r_k$ defined in \eqref{eq:ipm_ls}, but use the index~$\star$ instead of~$k$ to indicate that all functions are evaluated at $w^\star$, similarly for $\mathcal{M}_\star$.
We call solving~\eqref{eq:forward_ipm_sens} forward sensitivity computation.
Note that this linear system has the same structure as~\eqref{eq:ipm_ls} and can be solved using a factorization of the reduced matrix $\widetilde{\mathcal{M}}_\star$ in \eqref{eq:m_tilde}.
It is possible to obtain the sensitivities only for a subvector of $\theta$ by performing the backsolve~\eqref{eq:forward_ipm_sens} for the corresponding columns of $J_\star$.

For a given adjoint seed $\nu \in \R^{\npd}$, the adjoint sensitivity of the solution map is
\begin{align}
	\sadj\T & \coloneqq \nu\T \tfrac{\partial w^\sol_\textsc{ipm}}{\partial \theta}(w^\star;\tau, \theta) = -\nu\T \mathcal{M}_\star(w^\star; \tau, \theta)\inv J_\star(w^\star; \tau, \theta).
	\label{eq:adj_via_forw}
\end{align}
Transposing both sides yields
\begin{align}
	\label{eq:sadj_trans}
	\sadj & = -J_\star(w^\star; \tau, \theta)\T (\mathcal{M}_\star(w^\star; \tau, \theta)^{\shortminus \top} \nu).
\end{align}
Notice that a system of linear equations similar to \eqref{eq:ipm_ls} but with the coefficient matrix $\mathcal{M}_\star\T$ has to be solved.
Although $\mathcal{M}_\star$ is not symmetric, once $\mu$ and $s$ are eliminated, a system with the same lower dimensional symmetric coefficient matrix $\widetilde{\mathcal{M}}_\star$ is obtained.
This reduced system can be tackled with the same efficient algorithm used for the QP solution and the forward sensitivities, while the block elimination of ${\mathcal{M}}_\star\T$ requires a slightly different algorithm to process the right-hand side.

In the context of backpropagation through a neural network in which one layer is the solution of an optimization problem, the seed vector $\nu$ corresponds to the activations of the next layer. %

\begin{remark}[Approximate Hessian pitfall]
\label{rem:ift_pitfall}
At first glance the factorization of the matrix $\mathcal{M}_\star$, or to be more precise the factorization of $\widetilde{\mathcal{M}}_\star$, is naturally available at the solution when employing an interior-point NLP solver or an SQP solver that tackles the QP subproblems with an interior-point QP solver as outlined in Section~\ref{sec:sqp-ipm}.
An efficient implementation thus only requires the evaluation of $J_\star$ and a backsolve routine to complete the solution sensitivity computation, if these conditions hold.
However, the exact Hessian factorization is rarely available in practical NLP solvers.
Firstly, it is common to use a Hessian approximation, see Remark~\ref{remark:hess_approx}.
Secondly, the exact Hessian or its approximation are often regularized, e.g. by adding a positive constant onto the diagonal, commonly known as adding a Levenberg-Marquardt term \citep{Nocedal2006},
or by applying more advanced regularization techniques \citep{Verschueren2017}.

We emphasize the requirement of using an exact Hessian in the QP subproblem in the sensitivity computation referring to Thm.~\ref{thm:nlp-qp}.
If this assumption is not satisfied and a Hessian approximation is used instead, as suggested in previous works~\citep{Amos2018, Adabag2026}, the solution sensitivities might be highly degraded, as shown in Example~\ref{example:solution_sens_pendulum}.
\end{remark}

\textbf{Computational complexity.}\label{sec:forw_vs_adj}
Both the forward and adjoint sensitivity computations require computing $J_\star$ and $\mathcal{M}_\star$, as well as factorizing $\widetilde{\mathcal{M}}_\star$.
It is possible to recover the full solution sensitivity %
by solving~\eqref{eq:sadj_trans} using all $n_w$ canonical unit vectors for~$\nu$.
Compared to evaluating \eqref{eq:forward_ipm_sens} directly, this comes at the additional expense of a multiplication with $J_\star\T$.

For optimization problems with many parameters, i.e. $n_\theta \geq \npd $ and if only one or a few adjoint sensitivities are needed, it is computationally beneficial to compute these adjoints via \eqref{eq:sadj_trans} and to avoid computing the full forward sensitivity matrix in \eqref{eq:forward_ipm_sens}. %
In particular, this is most efficient when only one adjoint sensitivity is needed as is the case when performing backpropagation through a neural network in which one layer is the solution of an optimization problem.

\section{Differentiable MPC}\label{sec:diff_mpc}

In this section, we discuss the efficient solution of NLPs with an optimal control problem (OCP) structure using the tools introduced in the previous section.
In nonlinear MPC, one solves at each time instant an OCP of the form
\vspace{-2mm}
\begin{mini!}
	{\substack{x_0,\ldots, x_N, \\ u_0,\ldots, u_{N\!\shortminus 1} }}
	{\sum_{n=0}^{N\shortminus 1} L_n(x_n, u_n; \theta)  + M(x_N; \theta)}
	{\label{eq:acados_OCP}}
	{}
	\addConstraint{x_0}{= \bar x_0}
	\addConstraint{\!\!\!\!\!\!\!\!\!\!\!\!x_{n+1}\label{eq:acados_OCP_eq}}{=\phi_n(x_n, u_n; \theta),}{n=0,\ldots,N\!\shortminus 1}
	\addConstraint{0}{\geq h_n(x_n, u_n; \theta), \label{eq:acados_OCP_ineq}}{n=0,\ldots,N\!\shortminus 1}
	\addConstraint{0}{\geq h_N(x_N; \theta)},
\end{mini!}
where $\bar x_0\in\R^{\nx}$ is the initial state, $x_n\in \R^{\nx}$, $n=0,\dots, N$, represent a discrete state trajectory, and $u_n\in\R^{\nctrl}$, $n=0,\dots, N\shortminus1$, are control inputs.
The dynamics $\phi_n$ describe evolution of the state from $x_n$ to $x_{n+1}$.
The functions $h_n$ gather the constraints on $x_n$ and $u_n$, $n=0,\dots, N\shortminus1$, and $h_N$ represents the constraints on the terminal state.
The parameter vector $\theta\in \R^{n_\theta}$ represents all parameters that can change and with respect to which we want to compute the solution sensitivities.
Note that the OCP-structure exploiting algorithms and their implementations described in the sequel are applicable to problems with stage-wise varying dimensions, as described in \citet{Frey2025}.

\textbf{Structure exploiting solvers.}
While OCP \eqref{eq:acados_OCP} can be phrased as a general dense NLP as in~\eqref{eq:z_star_def}, exploiting the specific OCP structure in a tailored algorithm typically is most efficient for their numerical solution~\citep{Kouzoupis2018}.
When applying an interior-point QP solver to the OCP-structured problem~\eqref{eq:acados_OCP}, the Riccati-factorization can be applied to efficiently solve the linear system with coefficient matrix $\widetilde{\mathcal{M}}$ \citep{Steinbach1994, Frison2015a}.
In addition to the problem structure in \eqref{eq:acados_OCP}, some solvers, such as \acados~and \texttt{HPIPM}, offer special handling of slack variables.
These variables are structurally similar to control inputs, but do not enter the dynamics and can only enter the cost linearly or quadratically with diagonal Hessian. %
This structure often occurs in practice, especially when formulating soft constraints.

\textbf{Efficient solution and sensitivity computation.}
\label{sec:sol_sens_alg}
The computations required to compute the solution sensitivities of an OCP-structured NLP~\eqref{eq:z_star_def} can be summarized as follows:
\begin{enumerate}[nolistsep]
	\item[(S1)] Solve \eqref{eq:z_star_def} or a smoothed variant of its KKT conditions with $\tau\ind{min}$ to obtain $w^\sol_{\textsc{IPM}}(\bar{\theta})$. \label{step:nlp_solve}
	\item[(S2)] Evaluate the coefficient matrix $\mathcal{M}_\star$ with exact Hessian $Q_\star = \nabla^2_{z z}\mathcal{L}$. \label{step:eval_ex_hess}
	\item[(S3)] Factorize the corresponding matrix $\widetilde{\mathcal{M}}_\star$.
\end{enumerate}
The above steps need to be carried out independent of whether forward or adjoint sensitivities are computed.
For forward sensitivities, the remaining steps are:
\begin{enumerate}[nolistsep]
	\item[(S4f)] Evaluate $J_\star$, the derivatives of the residual map with respect to the parameters. \label{step:eval_J}
	\item[(S5f)] Solve the linear system \eqref{eq:forward_ipm_sens} to obtain the forward sensitivity result. \label{step:sens_solve}
\end{enumerate}

For adjoint sensitivities, the remaining steps are:
\begin{enumerate}[nolistsep]
	\item[(S4a)] Solve the linear system in brackets of \eqref{eq:sadj_trans}. \label{step:sens_solve_adj}
	\item[(S5a)] Compute the adjoint sensitivities with the seed computed in (S4a), without computing $J_\star$. \label{step:eval_J_adj}
\end{enumerate}

\section{Efficient implementation in \texttt{acados}}\label{sec:acados_sol_sens}
In this section, we present the efficient open-source implementation of the algorithms conceptually described above within the \acados{} software package.

\textbf{Two-solver approach.}
In order to alleviate the potential pitfalls pointed out in Sec.~\ref{sec:dense_solution_sens}, we suggest a two-solver approach to carry out the steps outlined in Sec.~\ref{sec:sol_sens_alg}.
This approach consists of
\begin{enumerate}[nolistsep] %
	\item A \emph{nominal solver}, which uses the most suitable Hessian approximation and regularization technique to obtain $w^\star(\bar{\theta})$ by performing step~(S1).
	\item A \emph{sensitivity solver}, which carries out steps (S2) to (S5a)/(S5f) of the algorithm above.
\end{enumerate}
The nominal solver could use any QP solver, cf. Remark~\ref{remark:riccati}, or even an optimization solver outside \texttt{acados}, as only the primal-dual solution~$w^\star(\bar{\theta})$ needs to be transferred to the sensitivity solver.
If choosing $\tau\ind{min}$ larger than the desired tolerance, the nominal solver needs to compute the primal-dual solution of the smoothed KKT conditions which can be obtained by an SQP variant with a suitable IPM QP solver, like \texttt{HPIPM}, or an NLP IPM solver like \texttt{IPOPT}~\citep{Waechter2006}.

Even if all steps are carried out by \acados, the flexibility of this approach allows for an overall more efficient algorithm, as structurally different techniques can be used for the two solvers.
Firstly, regularization techniques and Hessian approximations can be beneficial for the solution, as motivated in Sec.~\ref{sec:dense_solution_sens}, which can require explicit definitions of convex-over-nonlinear structures.
Secondly, it might be attractive to choose different QP solvers for the subproblems, including the (partial) condensing algorithm~\citep{Axehill2015, Frison2016}.
Condensing is especially favorable if multiple QP solver iterations (each requiring a KKT matrix factorization within an IPM) are performed for a single condensed QP, as typically in~(S1), but not in~(S5a)/(S5f).
In Appendix~\ref{sec:code_snippet}, we provide a code snippet showing how these steps can be performed in the proposed implementation.

\textbf{Efficient linear system solves.}
The linear systems associated with the computation of the forward and adjoint sensitivities, as well as the computation of the Newton step inside the interior-point QP solver, can be reduced to the Riccati factorization of $\widetilde{\mathcal{M}}$, which is efficiently implemented in \texttt{HPIPM}~\citep{Frison2020} utilizing the high-performance linear algebra library \texttt{BLASFEO}~\citep{Frison2018}.
For the computation of adjoint sensitivities, a routine was implemented in \texttt{HPIPM} to reduce the right-hand side of a linear system with coefficient matrix~$\mathcal{M}\T$ to one with~$\widetilde{\mathcal{M}}$.
The computational complexities of the Riccati factorization and the corresponding solution routine are $\mathcal{O}(N(\nx+\nctrl)^3)$ and $\mathcal{O}(N(\nx+\nctrl)^2)$, respectively.
\begin{remark}[Riccati variants \& regularity]
\label{remark:riccati}
Different variants of the Riccati factorization exist.
The classic Riccati recursion allows indefinite full-space Hessians and only requires positive definiteness of the reduced Hessian with respect to the dynamics~\eqref{eq:acados_OCP_eq}.
In addition, \texttt{HPIPM} provides a square-root implementation based on a Cholesky factorization, which has a lower associated computational cost, but requires the full-space Hessian to be positive definite~\citep{Frison2020a}.
Thus, the classic Riccati algorithm should be used for~(S5a)/(S5f), while the square-root algorithm can be leveraged in an efficient implementation of~(S1).
\end{remark}

\textbf{Parallelization.}
In order to utilize many CPU cores when computing the solution and sensitivities for a series of values $\bar{\theta}$, we implemented a batch solver class in \texttt{Python} based on a \texttt{C} implementation with \texttt{OpenMP} parallelization.
This class also implements parallelized solver interactions, such as updating parameters, and storing and loading initializations.

\textbf{Limitations.}
The proposed implementation is currently limited to CPU-based computation, reflecting the design philosophy of \texttt{acados}, which is implemented in \texttt{C} with minimal dependencies to support deployment on embedded platforms.
While GPU-based implementations may be attractive for training networks that include optimization layers, efficient CPU execution remains essential for real-time MPC on embedded controllers, where GPUs are uncommon due to cost and energy constraints and are often inefficient for small-scale linear algebra.
Moreover, nonlinear and constrained formulations hinder batch parallelization due to high variance in solution times, which partly explains the lack of GPU support (cf. Table~\ref{tab:diff_mpc}).
Another limitation is that the system dynamics~$\phi_n$ must be provided in discrete time using \texttt{CasADi} expressions~\citep{Andersson2019}.
Support for general parametric solution sensitivities within the \texttt{acados} integrators, which are used to efficiently solve initial value problems and can be embedded within OCP solvers~\citep{Frey2019,Frey2023}, is left for future work.

\textbf{Integration with ML frameworks.}
Our implementation is wrapped for use in common ML frameworks, with additional support for handling diverse parameters.
These functionalities are available in the \texttt{leap-c} project~\citep{leap-c_releases}; see~\ref{app:leap_c} for details.
An example of learning-enhanced MPC using this pipeline is presented in~\ref{sec:cartpole}.
End-to-end learning-enhanced MPC experiments that leverage the proposed differentiable nonlinear MPC implementation are described in \citet{Llanes2026, Jahncke2026}.

\section{Numerical examples}
\label{sec:examples}
This section presents three numerical examples.
The first example showcases the proposed implementation's ability to differentiate the solution of an OCP with parametric cost, dynamics and constraints.
In a second example, we compare the efficiency of the proposed SBQP solution sensitivity computation with state-of-the-art methods on a linear-quadratic OCP with simple input bounds.
Thirdly, a nonlinear OCP with high dimensional parameter vector allows us to evaluate the computation times of computing forward and adjoint solution sensitivities in \texttt{acados}.
All code required to reproduce the results in this paper is available, see Appendix~\ref{sec:reproducibility}.

\begin{example}[OCP with highly nonlinear parameter dependency]
\label{example:solution_sens_pendulum}
\begin{figure*}
	\centering
	  \includegraphics[width=0.99\linewidth]{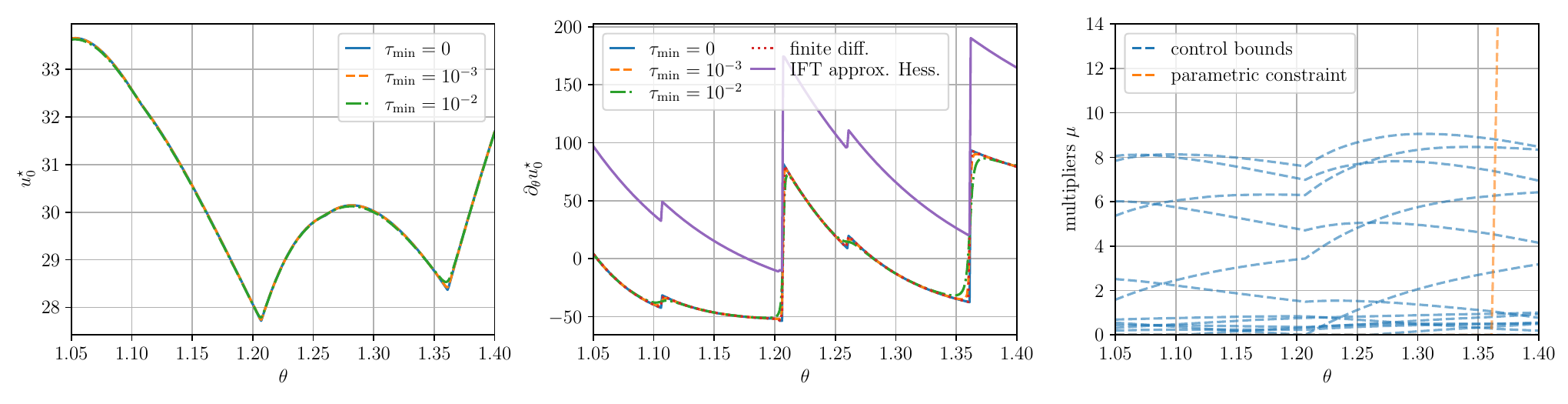}
	  \vspace{-2mm}
	  \caption{Optimal solution $u_0$ and sensitivities obtained by \acados~compared with finite differences for the OCP described in Example~\ref{example:solution_sens_pendulum}.
	  At active set changes, which are indicated by multipliers switching from being zero to being positive (or vice versa), the solution map is nondifferentiable.}
	  \label{fig:solution_sens_pendulum}
	  \vspace{-2mm}
\end{figure*}
We consider a nonlinear OCP with $\nx\!=\!4, \nctrl\!=\!1, n_\theta \!=\! 1$ and $N\!=\!50$.
The parameter $\theta$ enters the cost, the inequalities linearly and the dynamics nonlinearly.
A detailed description of the example is provided in Appendix~\ref{app:highly_param_ocp}.
Figure~\ref{fig:solution_sens_pendulum} visualizes the optimal solution $u_0^\star(\theta)$, as well as the solution sensitivities obtained via \acados~with different values of $\tau\ind{min}$.
The solver tolerance is set to $10^{\shortminus8}$.
Thus, for $\tau\ind{min} = 0$, the interior-point KKT system is solved with this accuracy, while the barrier parameter values $\tau_k$ never reach zero.
Additionally, sensitivities obtained by applying finite differences and sensitivities computed with the IFT approach and a Gauss-Newton Hessian approximation are plotted.
The solution sensitivities obtained with $\tau\ind{min} = 0$ are consistent with the ones obtained by finite differences, while the results obtained by the IFT approach with a Gauss-Newton Hessian approximation are completely off.
Additionally, the smoothing effect of larger values for $\tau\ind{min}$ is visible for both the solution and the corresponding sensitivities.
The kinks in the solution map are very pronounced and align with the active set changes indicated by multipliers changing from zero to a positive value. %
\end{example}

\begin{example}[{Benchmark example}]
\label{example:mpc_pytorch_benchmark}
{
\footnotesize
\begin{table*}
\centering
\caption{Timings in [ms] for solving $n_{\mathrm{batch}} \!=\! 128$ bounded LQR problems with varying parameters. In the config column, the difference to the default configuration is shown. In parentheses are multiples of the \acados{} runtime.\label{tab:mcp_pytorch}
}
\vspace{-1.5mm}
\footnotesize
\begin{tabular}{ccccccc}
\toprule
& \multicolumn{3}{c}{\textbf{Nominal solution}} & \multicolumn{3}{c}{\textbf{Solution + adjoint sens.}}\\
config & \acados & \texttt{mpc.pytorch} & \texttt{cvxpygen}& \acados & \texttt{mpc.pytorch} & \texttt{cvxpygen}\\
\midrule
default& 8.6 & $59 \;\,$ ($\times6.9$) & $235 \;\,$ ($\times27$)& 21.0 & $102 \;\,$ ($\times4.9$) & $329 \;\,$ ($\times16$)\\
$u_{\max} = 1$& 15.7 & $1697 \;\,$ ($\times110$) & $1140 \;\,$ ($\times73$)& 23.4 & $1793 \;\,$ ($\times77$) & $1239 \;\,$ ($\times53$)\\
$N = 100$& 29.3 & $671 \;\,$ ($\times23$) & $925 \;\,$ ($\times32$)& 68.2 & $828 \;\,$ ($\times12$) & $2268 \;\,$ ($\times33$)\\
$n_x = 16$, $n_u = 8$& 19.7 & $89 \;\,$ ($\times4.5$) & $585 \;\,$ ($\times30$)& 40.7 & $153 \;\,$ ($\times3.8$) & $856 \;\,$ ($\times21$)\\
\bottomrule
\end{tabular}
\vspace{-2mm}
\end{table*}

}

Next, we compare the proposed differentiable OCP solver with \texttt{mpc.pytorch}~\citep{Amos2018} and \texttt{cvxpygen} paired with \texttt{OSQP}~\citep{Schaller2026,Stellato2020}.
To enable comparison with these methods, we consider an OCP with quadratic cost, linear discrete-time dynamics, and bounds on the control inputs, cf. Table~\ref{tab:diff_mpc}.
We treat the dynamics matrices and offset, as well as the cost matrices, as parameters~$\theta$.
The problem dimensions are chosen as $N \!=\! 20$, $n_x \!=\! 8$, $n_u\!=\! 4$ if not written otherwise.
We consider both a loose-bound case and a case with active control bounds.
A detailed formulation and motivations for design choices of the benchmark problem can be found in Appendix~\ref{sec:mpc_pytorch_gpu}.
We solve $n\ind{batch}$ problems with different initial states $\bar{x}_0$, whose components are sampled from a standard normal distribution.
Table~\ref{tab:mcp_pytorch} reports nominal solution times and joint solution-plus-adjoint-sensitivity times.
Since the OCPs are strictly convex, it is not necessary to use the two-solver approach in \acados.
The table reports more than three times lower runtimes for \acados{} in every configuration compared to the other solvers.
We verified that all solvers converge to the same primal solution.
With tight control bounds, \texttt{mpc.pytorch} shows a large increase in runtime while the \texttt{acados} runtime changes much less.
Thus \texttt{acados} achieves speedup factors greater than 50 in this case.
The online computations needed when applying a (learned) MPC scheme correspond to the nominal solution, such that those timings can be seen as the inference time. %
The case $N=100$ shows runtime growth for \acados{} that is roughly consistent with linear scaling in the horizon.
The larger-dimension case shows a qualitatively similar runtime comparison to the default case.
Appendix~\ref{sec:mpc_pytorch_gpu} also reports a GPU-enabled variant of the benchmark. The loose-bound results show a similar qualitative comparison; for tight bounds, however, the GPU results use an earlier dynamics setup for which \texttt{mpc.pytorch} fails to converge, so the reported runtime ratios do not represent speedups for successful solves.
The main takeaway is that \acados~can handle more general problem formulations than the investigated alternatives while reducing the computation time.
\end{example}

\begin{figure}
	\centering
	\includegraphics[width=\columnwidth]{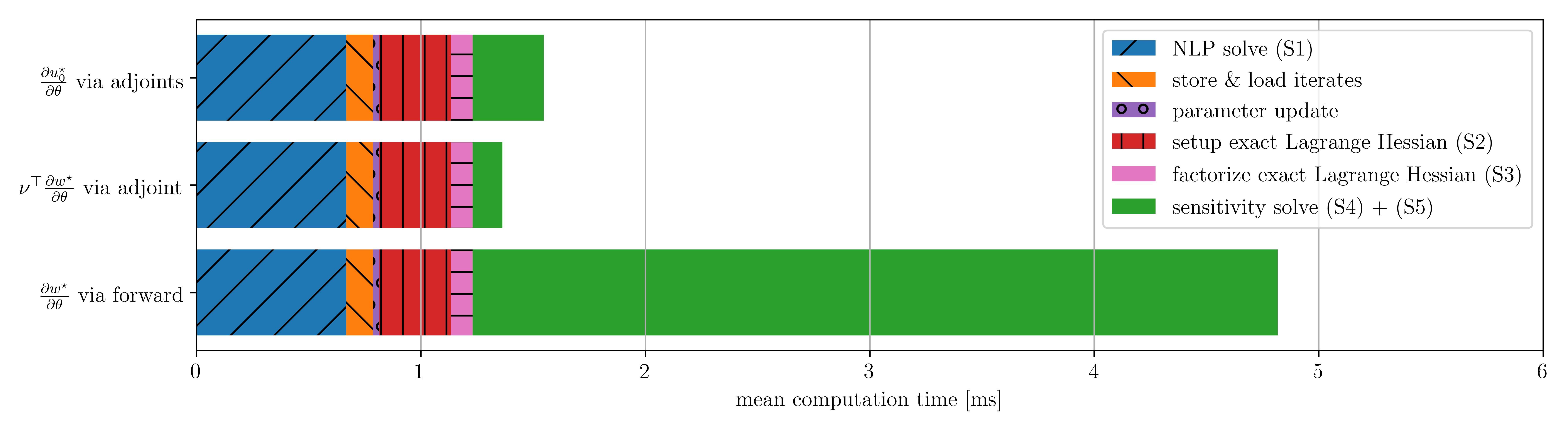}
	\vspace{-.4cm}
	\caption{Computation time for evaluating different kinds of sensitivities for Example~\ref{example:chain_sens}.}
	\vspace{-.5cm}
	\label{fig:chain_fwd_adj_timings}
\end{figure}

\begin{example}[Chain of masses]
\label{example:chain_sens}
We regard the example of a nonlinear OCP with a chain of masses model and many parameters, as described in \citet{Reinhardt2024a} with prediction horizon $N \!=\! 40$ and number of masses fixed to $n\ind{mass}\!=\!3$ such that $\nx \!=\! 9, \nctrl \!=\! 3, n_\theta \!=\! 113$.
The findings in~\citet{Reinhardt2024a} show a speedup of over one magnitude in forward sensitivity computation when using \texttt{acados}{} over general purpose libraries.
We consider different variants of solution sensitivity computation and analyze their computation times in Figure~\ref{fig:chain_fwd_adj_timings}.
We observe that the evaluation of one or a few adjoint sensitivities is computationally much cheaper compared to a full forward Jacobian sensitivity.
However, if many directions are needed, the computational cost of computing many adjoint direction of $J_\star$, see \eqref{eq:sadj_trans} and step (S5a), makes a computation via adjoints more expensive.
All variants only differ in the last step, namely the sensitivity solve with a given factorization.
The takeaway is that computing a single adjoint sensitivity is roughly 4 times faster than evaluating the full forward sensitivity and the adjoint sensitivity contains all required information when performing backpropagation to optimize parameters in a learning framework.

\end{example}

\section{Conclusion and outlook}\label{sec:conclusion}
This paper presents an efficient, modular, open-source implementation for solving general optimal control structured nonlinear programs and computing the corresponding forward or adjoint solution sensitivities within the software framework \texttt{acados}.
Various examples highlight its usability, flexibility and superior performance with respect to existing implementations.
This work builds the basis for integrating efficient and reliable solvers for optimal control problems as differentiable layers in machine learning networks.
Future work includes a comparison of the smoothed sensitivities with conservative Jacobians \citep{Bolte2021} with respect to their performance within a learning pipeline, as well as extending the \acados~implementation to allow using its efficient integrators and implementing additional regularizations to obtain sensitivities of locally non-unique solutions as suggested in \citet{Zuliani2025a}.

\anonymousOrSubmission{}
{
\subsection*{Acknowledgements}
This research was supported by DFG via projects 504452366 (SPP 2364), 560056112 (robust MPC), 535860958 (ALeSCo) and 525018088 (MAWERO).
}

\bibliographystyle{icml2026}
\bibliography{syscop}

\clearpage

\newpage
\clearpage
\appendix

\section{Appendix}

\newcommand{\yes}{\includegraphics[width=10pt]{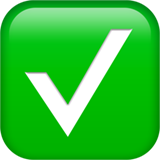}}
\newcommand{\no}{\includegraphics[width=10pt]{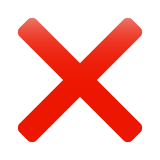}}
\newcommand{\neutral}{\includegraphics[width=10pt]{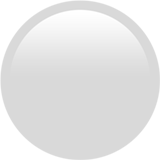}}
\newcommand{\unclear}{\includegraphics[width=10pt]{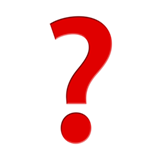}}
\newcommand{\tabcell}[1]{\begin{tabular}[c]{@{}c@{}} #1 \end{tabular}}
\newcommand{\ltabcell}[1]{\begin{tabular}[c]{@{}l@{}} #1 \end{tabular}}
\newcommand{\negspace}{\hspace{-.2cm}}
\begin{table*}
\centering
\begin{threeparttable}
\caption{
	Overview on differentiable optimization software and suitability for optimal control.
	\emph{OCP tailored} indicates whether the internal algorithms exploit the OCP sparsity structure.
	The column \emph{adjoint} indicates if adjoint sensitivities as needed for backpropagation are implemented.
	The column \emph{ineq. support} indicates whether inequalities are supported.
	The following columns summarize which parts of the problem can be nonlinear.
\label{tab:diff_mpc}
}
\footnotesize
\begin{tabular}{ccccccccccl}
\toprule
Software & OCP & ad- & ineq. & \multicolumn{3}{c}{Nonlinearity in} & \tabcell{exact} & parallel- & \multicolumn{2}{l}{permissive} \\
/ publication & tailored & joint & support & cost & eq. & ineq. & \tabcell{results} & ization & \multicolumn{2}{l}{license} \\
\midrule
\tabcell{\texttt{acados} \\this paper} & \tabcell{\yes} & \tabcell{\yes} & \tabcell{\yes} & \tabcell{\yes} & \tabcell{\yes} & \tabcell{\yes} & \tabcell{\yes} & \tabcell{\texttt{OpenMP}} & \tabcell{\yes} & \negspace BSD-2 \\ \addlinespace[5pt]
\tabcell{\texttt{mpc.pytorch} \\ \citet{Amos2018}} & \tabcell{\yes} & \tabcell{\yes} & \tabcell{\no\tnote{a}}  & \tabcell{\yes} & \tabcell{\no} & \tabcell{\no} & \tabcell{\no} & \texttt{PyTorch} & \tabcell{\yes} & \negspace MIT \\ \addlinespace[5pt]
\tabcell{\texttt{cvxpygen +} \\ \texttt{cvxpylayers}\\ \citet{Schaller2022}\\ \citet{Schaller2026}} & \tabcell{\no} & \tabcell{\yes} & \tabcell{\yes} & \tabcell{\neutral \tnote{b}} &  \tabcell{\no} & \tabcell{\neutral \tnote{b}}  & \tabcell{\yes} & \tabcell{\texttt{Python}\\ multi-\\threading} & \tabcell{\yes} & \negspace \ltabcell{Apache\\ 2.0} \\ \addlinespace[5pt]
\tabcell{\texttt{diffmpc} \\(Toyota)\\ \citet{Adabag2026}} & \tabcell{\yes} & \tabcell{\yes} & \tabcell{\no} & \tabcell{\yes} & \tabcell{\yes} & \tabcell{\no} & \tabcell{\no} & \texttt{JAX} & \tabcell{\no} & \negspace \ltabcell{CC BY-\\NC-ND\\ 4.0} \\ \addlinespace[5pt]
\tabcell{\texttt{trajax} \\ (Google) \\ no paper}& \tabcell{\yes} & \tabcell{\yes} & \tabcell{\no} & \tabcell{\yes} & \tabcell{\yes} & \tabcell{\no}  & \neutral \tnote{c} & \texttt{JAX} & \tabcell{\yes} & \negspace \ltabcell{Apache\\ 2.0} \\ \addlinespace[5pt] %
\tabcell{\texttt{theseus} \\ (Facebook) \\ \citet{Pineda2022}} & \tabcell{\no} & \tabcell{\yes} & \tabcell{\no} & \tabcell{\neutral \tnote{d}} & \tabcell{\no} & \tabcell{\no} & \tabcell{\no} & \texttt{PyTorch} & \tabcell{\yes} & \negspace MIT \\ \addlinespace[5pt]
\tabcell{\texttt{sIPOPT} \\ \citet{Pirnay2012}\tnote{e}} & \tabcell{\no} & \tabcell{\no} & \tabcell{\yes} & \tabcell{\yes} & \tabcell{\yes} & \tabcell{\yes} & \tabcell{\yes} & -- & \tabcell{\yes} & \negspace EPL \\ \addlinespace[5pt]
\tabcell{\texttt{CasADi} \\ native \\ \citet{Andersson2018a}} & \tabcell{\no} & \tabcell{\yes} & \tabcell{\yes} & \tabcell{\yes} & \tabcell{\yes} & \tabcell{\yes} & \tabcell{\yes} & -- & \tabcell{\yes} & \negspace LGPL \\ \addlinespace[5pt]
\bottomrule
\end{tabular}

\begin{tablenotes}
\footnotesize
\item[a] while simple control bounds can be added to the problem formulation, the experiments in Example~\ref{example:mpc_pytorch_benchmark} and Appendix~\ref{sec:mpc_pytorch_gpu} revealed that the solver fails to converge if active set is nontrivial.
\item[b] only convex cost and constraints.
\item[c] A set of options appears to allow the use of the exact Hessian for solution sensitivity computation; however, the software and documentation do not specify the conditions under which the resulting sensitivities are accurate.
\item[d] supports only nonlinear least-squares cost and no equality constraints, as needed for dynamics. OCPs with quadratic costs can be implemented using the single-shooting formulation.
\item[e] functionality extremely limited in practice. Not usable in \texttt{CasADi} and abandoned rewrite effort in \url{https://github.com/athrpf/sipopt2}.
\end{tablenotes}
\end{threeparttable}
\end{table*}

\subsection{Differentiable MPC software}
\label{sec:diff_mpc_software}
Table~\ref{tab:diff_mpc} provides a qualitative comparison of differentiable optimization solvers, with a focus on optimal control.
It can be seen that inequality treatment, a feature that is essential to leverage safety benefits of learning-enhanced MPC, is scarce.
In particular, no OCP-tailored software other than the proposed implementation supports inequality constraints.
Some existing solutions rely on identifying the active set at the solution and applying the IFT, e.g. the ``native'' \texttt{CasADi} solution sensitivity implementation~\citep{Andersson2018a}\footnote{\url{https://github.com/casadi/casadi/discussions/3156}}.
This clear separation does not allow one to differentiate the solution across active set changes in contrast to the proposed approach.

\subsection{Implicit function theorem} \label{sec:ift}
As it is essential for this work, we state the implicit function theorem (IFT) in its classical form~\citet{Dini1907}.
\begin{theorem}[Implicit function theorem (IFT)] \label{thm:ift}
Let $r: \R^{n_\theta} \times \R^{n_w} \to \R^{n_w}$ be continuously differentiable in a neighborhood of $(\bar{\theta}, \bar{w})$ with $r(\bar{\theta}, \bar{w}) = 0$, and let the Jacobian $\nabla_w r(\bar{\theta}, \bar{w})$ be nonsingular.
Then the solution mapping
\begin{align}
	S: \theta \mapsto \{ w \in \R^{n_w} \, | \, r(\theta, w) = 0 \}
\end{align}
has a single-valued localization $s$ around $\bar{\theta}$ for $\bar{w}$ which is continuously differentiable in a neighborhood $Q$ of $\bar{\theta}$ with Jacobian satisfying
\begin{align}
	\nabla s (\theta) = - \nabla_w r(\theta, s(\theta))\inv \nabla_\theta r(\theta, s(\theta)) \; \textrm{for all} \; \theta\in Q.
\end{align}
\end{theorem}

\subsection{Interpreting SBQP as Sequential Convex Programming (SCP) applied to the barrier problem}
\label{sec:ipm_sqp_details}
This section presents additional details on SBQP introduced in Section~\ref{sec:sqp-ipm}, in particular if a value $\tau\ind{min} > 0$ is used.
For $\tau\ind{min} = 0$, the algorithm can be interpreted as a standard SQP method.

In the following, we assume $\tau \coloneqq \tau\ind{min} > 0$ is used and show that SBQP solves the nonlinear barrier problem~\eqref{eq:dense_nlp_log_barrier}.
Moreover, we show that, in some special cases, it can be interpreted in the framework of sequential convex programming (SCP).

Let us repeat the original NLP
\begin{mini}
	{\substack{z\in\R^{n_z}}}
	{f(z; \theta)}
	{\label{eq:dense_nlp_appendix}}
	{}
	\addConstraint{g(z; \theta)}{= 0}
	\addConstraint{h(z; \theta)}{\leq 0.}
\end{mini}
Replacing the inequality constraint by a barrier term, introducing the barrier parameter $\tau > 0$, yields the following log-barrier problem:
\begin{mini}
	{\substack{z\in\R^{n_z}}}
	{f(z; \theta) - \tau\sum_{i=1}^{n_h} \log(-h_i(z;\theta))}
	{\label{eq:dense_nlp_log_barrier_no_slack}}
	{}
	\addConstraint{g(z; \theta)}{= 0,}
\end{mini}
with associated KKT conditions
\begin{subequations}
\label{eq:nlp_kkt_ipm}
\begin{align}
	\grad{z}  f(z; \theta) + \grad{z} g(z; \theta) \lambda &+ \grad{z} h(z; \theta) \mu = 0, \\
	g(z; \theta) &= 0, \\
	h(z; \theta) &\leq 0,\\
	\mu & \geq 0, \\
	\mu_i h_i(z; \theta) &= -\tau\ind{min}, \, i=1,\dots,n_h.
\end{align}
\end{subequations}
Their derivation is similar to the ones derived next with linearized constraints and a quadratic cost approximation.
We also call \eqref{eq:nlp_kkt_ipm} the smoothed interior-point KKT system as introduced next to~\eqref{eq:KKT_comp_ipm}.

By introducing slack variables $s\in\R^{n_h}$, we can reformulate~\eqref{eq:dense_nlp_log_barrier_no_slack} as
\begin{mini}
	{\substack{z\in\R^{n_z}, s\in\R^{n_h}}}
	{f(z; \theta) - \tau\sum_{i=1}^{n_h} \log(s_i)}
	{\label{eq:dense_nlp_log_barrier}}
	{}
	\addConstraint{g(z; \theta)}{= 0}
	\addConstraint{h(z; \theta) + s}{= 0.}
\end{mini}

Let us only keep the convexity of the logarithmic barrier term, quadratically approximate the cost function and linearize the remaining parts of the problem~\eqref{eq:dense_nlp_log_barrier} to arrive at
\begin{mini}
	{\substack{\Delta z\in\R^{n_z}, \\ s\in\R^{n_h}}}
	{\!\!\! f_k \!+\! q_k\T \Delta z \!+\! \tfrac{1}{2} \Delta z\T \nabla^2_{zz} f(z_k; \theta) \Delta z \!-\! \tau\!\sum_{i=1}^{n_h} \log(s_i)}
	{\label{eq:scp_log_barrier}}
	{}
	\addConstraint{g_k + G_k \Delta z}{= 0}
	\addConstraint{h_k + H_k \Delta z + s}{= 0,}
\end{mini}
where the same shorthands as in Section~\ref{sec:sqp-ipm} are used, namely
$q_k = \grad{z} f(z_k; \theta)$, $g_k = g(z_k; \theta)$, $h_k = h(z_k; \theta)$ and $G_k = \nabla g(z_k; \theta)^\top $, $H_k = \nabla h(z_k; \theta)^\top $, as well as $f_k = f(z_k; \theta)$.

The KKT conditions of this problem read:
\begin{subequations}
	\label{eq:barrier_kkt}
	\begin{align}
	q_k + \nabla^2_{zz} f(z_k; \theta) \Delta z + G_k\T \lambda + H_k\T \mu &= 0, \\
	-\frac{\tau}{s_i} + \mu_i &= 0, \, i=1,\dots,n_h, \label{eq:barrier_kkt_stat_slack}\\
	g_k + G_k \Delta z &= 0, \\
	h_k + H_k \Delta z + s &= 0, \\
	s & > 0. \label{eq:barrier_kkt_slack_pos}
	\end{align}
\end{subequations}

Equation~\eqref{eq:barrier_kkt_stat_slack} corresponds to $\grad{s} \mathcal{L} = 0 $ and~\eqref{eq:barrier_kkt_slack_pos} is needed to ensure a well-defined barrier term.
We can rewrite~\eqref{eq:barrier_kkt_stat_slack} as $\mu_i s_i = \tau $ and deduce $\mu > 0 $.
This allows us to rewrite~\eqref{eq:barrier_kkt} as
\begin{subequations}
	\label{eq:barrier_kkt_rewrite}
	\begin{align}
	\nabla^2_{zz} f(z_k; \theta) \Delta z + q_k + G_k\T \lambda + H_k\T \mu &= 0, \label{eq:kkt_nlp_stat}\\
	g_k + G_k \Delta z &= 0, \\
	h_k + H_k \Delta z + s &= 0, \\
	\mu_i s_i &= \tau, \, i=1,\dots,n_h, \label{eq:kkt_nlp_comp}\\ 
	s, \mu & > 0. %
	\end{align}
\end{subequations}
These conditions correspond to the ones solved in an IPM QP solver, namely \eqref{eq:ipm_stat}-\eqref{eq:ipm_comp}, with the only difference being that \eqref{eq:ipm_stat} uses a Hessian approximation of the Lagrangian $Q_k$, while \eqref{eq:kkt_nlp_stat} uses the exact Hessian of the objective.

These matrices coincide if the objective is of linear least-squares type and a Gauss-Newton Hessian approximation is used for $Q_k$.
For an NLP with least-squares objective, i.e.
$f(z; \theta) = \frac{1}{2} \|F(z; \theta) - y\ind{ref}\|_W^2$, the Gauss-Newton Hessian approximation is $Q_k=J_k\T W J_k $ with $J = \frac{\partial F}{\partial z}(z_k, \theta) $.
If the function $F$ is linear in $z$, the Gauss-Newton Hessian approximation corresponds to the exact Hessian of the objective $\nabla^2_{zz} f(z_k; \theta)$.

In case of a linear least-squares objective and a Gauss-Newton Hessian approximation, SBQP can be interpreted in the framework of Sequential Convex Programming (SCP)~\citep{Messerer2021a} on the nonlinear barrier problem~\eqref{eq:dense_nlp_log_barrier} and the linear convergence rate result, Thm 4.5 in \citet{Messerer2021a}, holds.

In the general case, where $f$ is not of linear least-squares type, and any Hessian approximation $Q_k$ is used in the QP solver, the algorithm still provides a solution of \eqref{eq:dense_nlp_log_barrier}, assuming it converges.

\subsection{Solution sensitivities for NLPs}\label{sec:solution_sens_nlp}

This section provides a proof of Theorem~\ref{thm:nlp-sens-exist}, which follows from Assumption~\ref{as:regularity} and Theorem~\ref{thm:ift}.

Strict complementarity allows us to isolate the strictly active inequality constraints at the solution $z^\star(\bar \theta)$ and denote them by $h_\mathcal{A}$, with associated multipliers $\mu_\mathcal{A}$.
All other inequalities are inactive with zero multiplier values.
We want to apply the IFT in Theorem~\ref{thm:ift} at the solution of the KKT conditions $(z^\star(\bar \theta), \lambda^\star(\bar{\theta}), \mu_{\mathcal{A}}^\star(\bar \theta))$ to the residual function
\begin{align} \label{eq:residuals}
	r(z, \lambda, \mu_\mathcal{A}; \theta) =
	\begin{bmatrix}
		\grad{z}  f(z; \theta) + \grad{z} g(z; \theta) \lambda + \grad{z} h_\mathcal{A}(z; \theta) \mu_\mathcal{A} \\
		g(z; \theta) \\
		h_\mathcal{A}(z; \theta)
	\end{bmatrix}\!.
\end{align}
The matrix to be inverted can be written as
\begin{align}
	&\frac{\partial r}{\partial (z, \lambda, \mu_{\mathcal{A}})} =
	\begin{bmatrix}
		\nabla^2_{zz} \mathcal{L}(z, \lambda, \mu; \theta) & \nabla_z g(z; \theta) & \nabla_z h_\mathcal{A}(z; \theta) \\
		\nabla_z g(z; \theta)\T & 0 & 0 \\
		\nabla_z h_\mathcal{A}(z; \theta)\T & 0 & 0
	\end{bmatrix}\!.
\end{align}
LICQ and SOSC imply that this matrix is invertible at the solution, cf. Lemma 16.1 in \citet{Nocedal2006}, which allows us to apply the IFT and concludes the proof of Theorem~\ref{thm:nlp-sens-exist}.

In particular, the IFT implies that the solution sensitivity can be computed as
\begin{align}
\begin{bmatrix}
	\tfrac{\partial z^\star}{\partial \theta} (\bar \theta)\\[6pt]
	\tfrac{\partial \lambda^\star}{\partial \theta} (\bar \theta) \\[6pt]
	\tfrac{\partial \mu_{\mathcal{A}}^\star}{\partial \theta} (\bar \theta)
\end{bmatrix}
= -
\left(\frac{\partial r }{\partial (z, \lambda, \mu_{\mathcal{A}})}\right)^{-1} \frac{\partial r}{\partial \theta},
\label{eq:solution_sens_fwd}
\end{align}
where the derivatives of $r$ are evaluated at the solution $(z^\star(\bar \theta), \lambda^\star(\bar{\theta}), \mu_{\mathcal{A}}^\star(\bar \theta))$.

\subsection{Proof: Barrier NLP and barrier QP sensitivities coincide}\label{sec:bqp_sens_match}
This section provides a proof of Theorem~\ref{thm:barrier-bqp}, showing that the optimality conditions of the barrier QP match the nonlinear smoothed KKT conditions at convergence of SBQP.

The linearized problem quantities of the BQP at the KKT point are defined as
$ q_\star = \grad{z} f(z^\star; \theta)$, $Q_\star = \nabla^2_{zz}\mathcal{L}(z^\star, \lambda^\star, \mu^\star; \bar\theta)$,
$ g_\star = g(z^\star; \theta)$, $ h_\star = h(z^\star; \theta)$, $ G_\star = \nabla g(z^\star; \theta)\T$, $ H_\star = \nabla h(z^\star; \theta)\T$.
The KKT conditions of~\eqref{eq:bqp} at the KKT point are
\begin{subequations}
	\label{eq:diffmpc_KKT_BQP}
	\begin{align}
		Q_\star \Delta z + q_\star + G_\star\T \lambda + H_\star\T \mu &= 0, \\
		g_\star + G_\star \Delta z &= 0, \\
		h_\star + H_\star \Delta z + s &= 0, \\
		\mu_i s_i &= \tau, \; i=1,\dots,n_h, \\
		s, \mu &> 0.
	\end{align}
\end{subequations}
Since $(z^\star, s^\star, \lambda^\star, \mu^\star)$ is a KKT point of~\eqref{eq:dense_nlp_log_barrier} at $\bar\theta$, the relations $g(z^\star; \bar\theta)=0$, $h(z^\star; \bar\theta)+s^\star = 0$, $\mu_i^\star s_i^\star = \tau$, and
\begin{align*}
	\nabla_{z z}^2 \mathcal{L}(z^\star, \lambda^\star, \mu^\star; \bar\theta) \cdot 0 + q_\star + G_\star\T \lambda^\star + H_\star\T \mu^\star = 0
\end{align*}
show that $(\Delta z, s, \lambda, \mu) = (0, s^\star, \lambda^\star, \mu^\star)$ satisfies~\eqref{eq:diffmpc_KKT_BQP}, which completes the first part.

The sensitivities of the BQP solution can be computed via the IFT in Theorem~\ref{thm:ift}, given LICQ and SOSC.
This yields
\begin{align}
	\setlength\arraycolsep{1pt}
	\frac{\partial}{\partial \theta}
	\begin{bmatrix}
		\Delta z^\sol_\mathrm{BQP} (\tau; \bar \theta, v^\star)\\
		\lambda^\sol_\mathrm{BQP} (\tau; \bar \theta, v^\star)\\
		\mu^\sol_\mathrm{BQP} (\tau; \bar \theta, v^\star)\\
		s^\sol_\mathrm{BQP} (\tau; \bar \theta, v^\star)
	\end{bmatrix}
	=
	-
	\begin{bmatrix}
		Q_\star & G_\star\T & H_\star\T & 0 \\
		G_\star & 0 & 0 & 0 \\
		H_\star & 0 & 0 & \eye \\
		0 & 0 & S^\star & M^\star
	\end{bmatrix}^{-1}
	r_{\mathrm{sens},\star}^{\mathrm{BQP}},
	\label{eq:solution_sens_BQP_linear}
\end{align}
where $S^\star = \mathrm{diag}(s^\star)$ and $M^\star = \mathrm{diag}(\mu^\star)$ and
\begin{align}
	r_{\mathrm{sens},\star}^{\mathrm{BQP}}
	=
	\frac{\partial}{\partial \theta}
	\begin{bmatrix}
		Q_\star \Delta z + q_\star + G_\star\T \lambda + H_\star\T \mu \\
		g_\star + G_\star \Delta z \\
		h_\star + H_\star \Delta z + s \\
		(\mu_i s_i - \tau)_{i=1,\dots,n_h}
	\end{bmatrix}.
\end{align}

This linear system coincides with \eqref{eq:forward_ipm_sens} if the exact Hessian $Q_\star = \nabla^2_{z z}\mathcal{L}(z^\star, \lambda^\star, \mu^\star; \bar\theta)$ is used, which concludes the proof.

\subsection{Optimality conditions of the quadratic subproblem at SQP convergence}\label{sec:qp_sens_match}
Since the proposed SBQP method corresponds to SQP when using $\tau\ind{min} = 0$, we also show under which conditions the NLP sensitivities match those of the QP subproblem at SQP convergence in the next theorem.

\begin{theorem}[NLP and QP sensitivities coincide] \label{thm:nlp-qp}
Suppose Ass.~\ref{as:regularity} holds at $z^\star$ with parameter $\bar{\theta}$, let $\lambda^\star, \mu^\star$ be the corresponding multipliers, and define $v^\star \coloneqq (z^\star, \lambda^\star, \mu^\star)$.
Regard the QP \eqref{eq:dense_QP} which is obtained at the solution $(z^\star, \lambda^\star, \mu^\star)$ 
and assume an exact Hessian
$Q_\star = \nabla^2_{zz}\mathcal{L}(z^\star, \lambda^\star, \mu^\star; \bar\theta)$ is used.
Then, the solution maps $z^\sol(\theta)$ and $z^\star + \Delta z_\mathrm{QP}^\sol(\theta, v^\star)$, and their sensitivities, $\frac{\partial z^\sol}{\partial \theta}(\theta)$ and $\frac{\partial \Delta z_\mathrm{QP}^\sol}{\partial \theta}(\theta, v^\star)$ coincide at $\bar{\theta}$.
The same holds for the solution maps of the Lagrange multipliers.
\end{theorem}
\begin{proof}
The KKT conditions of QP~\eqref{eq:dense_QP} read as follows:
\begin{subequations}
	\label{eq:KKT_QP}
	\begin{align}
		Q_k \Delta z\ind{QP} + q_k + G_k\T \lambda\ind{QP} + H_k\T \mu\ind{QP} &= 0, \\
		g_k + G_k \Delta z\ind{QP} &= 0, \\
		h_k + H_k \Delta z\ind{QP} &\leq 0, \\
		\mu_{\mathrm{QP},i} \cdot (h_k + H_k \Delta z\ind{QP})_i &= 0, \;\, i=1,\dots, n_h, \\
		\mu\ind{QP} &\geq 0.
	\end{align}
\end{subequations}
SOSC of~\eqref{eq:z_star_def} at $(z^\star, \lambda^\star, \mu^\star)$ implies SOSC of the QP~\eqref{eq:dense_QP} at $v\ind{QP}^\star = (\Delta z\ind{QP}, \lambda^\star, \mu^\star)$, with $\Delta z\ind{QP}=0$.
With the regularity properties in Ass.~\ref{as:regularity}, it follows that this is a unique global solution of QP~\eqref{eq:dense_QP}.
Thus, we have $ z^\star + \Delta z_\mathrm{QP}^\sol(\bar{\theta}, v^\star) = z^\star $.

Next, we analyze the optimality conditions of the QP and show that at SQP convergence, the solution sensitivities of the NLP coincide with the ones of the exact Hessian QP subproblem.
Assuming strict complementarity, the optimality conditions~\eqref{eq:KKT_QP} simplify to
\begin{subequations}
	\label{eq:KKT_QP_active}
	\begin{align}
		Q_k \Delta z_\mathrm{QP} + q_k + G_k\T \lambda\ind{QP} + H_{k, \mathcal{A}}\T \mu_{\mathrm{QP},\mathcal{A}} &= 0, \\
		g_k + G_k \Delta z_\mathrm{QP} &= 0, \\
		h_{k,\mathcal{A}} + H_{k,\mathcal{A}} \Delta z_\mathrm{QP} &= 0.
	\end{align}
\end{subequations}
The sensitivities of the QP solution $(\Delta z_\mathrm{QP}^\star, \lambda\ind{QP}^\star, \mu\ind{QP}^\star)$ can be computed via the IFT as in \ref{sec:solution_sens_nlp}, which yields
\begin{align}
	\setlength\arraycolsep{1pt}
	\begin{bmatrix}
		\tfrac{\partial \Delta z_\mathrm{QP}^\star}{\partial \theta} (\bar \theta)\\[6pt]
		\tfrac{\partial \lambda\ind{QP}^\star}{\partial \theta} (\bar \theta) \\[6pt]
		\tfrac{\partial \mu_{\mathrm{QP},\mathcal{A}}^\star}{\partial \theta} (\bar \theta)
	\end{bmatrix}
	= 
	-
	\begin{bmatrix} \setlength{\arraycolsep}{1pt}
		Q_k & G_k\T & H_{k, \mathcal{A}}\T \\
		G_k & 0 & 0 \\
		H_{k, \mathcal{A}} & 0 & 0
	\end{bmatrix}%
	^{-1}
	r_{\mathrm{sens},k},
	\label{eq:solution_sens_QP}
\end{align}
where
\begin{align}
r_{\mathrm{sens},k} = 	\frac{\partial}{\partial \theta}
\!\!
\begin{bmatrix}
	Q_k \Delta z_\mathrm{QP} + q_k + G_k\T \lambda\ind{QP} + H_{k, \mathcal{A}}\T \mu_{\mathrm{QP},\mathcal{A}} \\
	g_k + G_k \Delta z_\mathrm{QP} \\
	h_{k,\mathcal{A}} + H_{k,\mathcal{A}} \Delta z_\mathrm{QP}
\end{bmatrix}.
\end{align}
At the NLP solution, the active sets of the QP and the NLP coincide, and we have $H_k = \nabla h(z^\star;\bar{\theta})^\top, G_k = \nabla g(z^\star;\bar{\theta})^\top$.
If the QP uses an exact Hessian of the Lagrangian, i.e. $Q_k = \nabla^2_z \mathcal{L}(z_k, \lambda_k, \mu_k; \bar{\theta})$, it follows that the coefficient matrices in \eqref{eq:solution_sens_QP} and \eqref{eq:solution_sens_fwd} coincide at SQP convergence.
Since at SQP convergence, we have $\Delta z_\mathrm{QP} = 0$ and $g_k = g(z^\star; \bar{\theta}), h_k = h(z^\star; \bar{\theta})$ and $\lambda\ind{QP} = \lambda^\star, \mu\ind{QP} = \mu^\star$, also the right-hand side of the linear systems in \eqref{eq:solution_sens_QP} and \eqref{eq:solution_sens_fwd} coincide.

This shows, that at the NLP solution, the solution sensitivities of the NLP coincide with the one corresponding to the exact Hessian QP, which completes the proof of Theorem~\ref{thm:nlp-qp}.
\end{proof}

\subsection{Regularity assumptions of original and barrier NLP}
\label{app:reg_log_barrier}
In this section, we show that the proposed barrier smoothing approach allows us to weaken Assumption~\ref{as:regularity} and that assuming strong SOSC instead of SOSC plus strict complementarity is sufficient to obtain well-defined solution sensitivities of the barrier problem.
In particular, this makes the proposed method applicable to differentiate the smoothed solution map across active set changes of the original problem.

\begin{definition}[Strong second-order sufficient condition]
At a KKT point $(\bar z,\bar\lambda,\bar\mu)$ for parameter $\bar\theta$, let $\mathcal A_+\coloneqq\{i:\bar\mu_i>0\}$, $\bar G\coloneqq\nabla_z g(\bar z; \bar\theta)^\top$, and $\bar H\coloneqq\nabla_z h(\bar z; \bar\theta)^\top$.
We say strong SOSC holds at $(\bar z,\bar\lambda,\bar\mu)$ if
\begin{align*}
	d^\top\nabla^2_{zz}\mathcal L(\bar z,\bar\lambda,\bar\mu; \bar\theta)d>0
	\; \text{for all }d\in K_0\setminus\{0\},\qquad
	K_0\coloneqq\{d:\bar Gd=0,\ \bar H_i d=0\text{ for all }i\in\mathcal A_+\}.
\end{align*}
\end{definition}

\begin{theorem}\label{thm:barrier_differentiability}
Suppose at a solution $(\bar z, \bar \lambda, \bar \mu)$ to the original NLP ~\eqref{eq:z_star_def} for parameter value $\bar{\theta}$, LICQ and strong SOSC hold.
Crucially, strict complementarity might \emph{not} be satisfied.

Consider the barrier QP~\eqref{eq:bqp} formed at
$(\bar z,\bar\lambda,\bar\mu)$ and parameter $\bar\theta$, with
$q_k=\bar q\coloneqq\nabla_z f(\bar z;\bar\theta)$,
$Q_k=\nabla^2_{zz}\mathcal L(\bar z,\bar\lambda,\bar\mu;\bar\theta)$,
$g_k=\bar g$, $G_k=\bar G$, $h_k=\bar h$, $H_k=\bar H$, and
$\tau_{\min}=\tau>0$.
Assume its solution branch $(\Delta z^*(\tau),s^*(\tau),\lambda^*(\tau),\mu^*(\tau))$ converges to $(0,-\bar h,\bar\lambda,\bar\mu)$ as $\tau\downarrow0$. Note that Theorem~\ref{thm:smooth_KKT_sol} gives sufficient conditions for convergence of the nonlinear smoothed KKT branch.
Then, the solution map corresponding to \eqref{eq:bqp} is locally unique and differentiable with respect to $\theta$ for sufficiently small $\tau$.
\end{theorem}

\begin{proof}
Fix $\tau > 0$ on the assumed solution branch.
By Theorem~\ref{thm:nlp-sens-exist}, it suffices to show that the barrier QP satisfies LICQ and SOSC, as strict complementarity is trivially satisfied for the equality-constrained barrier QP.
LICQ of the original NLP implies that $\bar G$ has full row rank, and the identity block for the slack variables makes the barrier QP constraint Jacobian $\begin{bmatrix}\bar G & 0 \\ \bar H & I\end{bmatrix}$ full row rank.
It remains to establish SOSC for the barrier QP for sufficiently small $\tau$.

Let $W\coloneqq \nabla^2_{zz}\mathcal L(\bar z, \bar\lambda, \bar\mu; \theta)$ and $K\coloneqq \ker \bar G$.
The tangent directions $(d_z,d_s)$ of \eqref{eq:bqp} satisfy

$$
\bar G d_z = 0, \qquad \bar H d_z + d_s = 0.
$$

Hence each tangent direction is determined by $d \coloneqq d_z \in K$ with $d_s = -\bar H d$.
The Lagrangian Hessian is $\operatorname{diag}(W,\,\tau S^{-2})$, where $S\coloneqq \operatorname{diag}(s^*(\tau))$.
The corresponding reduced quadratic form is

\begin{align*}
	Q_\tau(d) \coloneqq d^\top W d + \tau \sum_{i=1}^{n_h}\frac{(\bar H_i d)^2}{s_i^*(\tau)^2}, \qquad d \in K.
\end{align*}

Thus SOSC holds exactly when $Q_\tau(d)>0$ for every nonzero $d\in K$.
By homogeneity, it suffices to check positivity on the unit sphere

\begin{align*}
\Sigma \coloneqq K \cap \{d : \|d\| = 1\}.
\end{align*}

Let $\mathcal A_+ \coloneqq \{i : \bar\mu_i > 0\}$.
Define the subspace $K_0$ and the function $\rho$ by

\begin{align*}
K_0 &\coloneqq \{d\in K : \bar H_i d = 0 \text{ for all } i\in\mathcal A_+\},\\
\rho(d) &\coloneqq \sum_{i\in\mathcal A_+}(\bar H_i d)^2, \qquad d\in K.
\end{align*}

Then $K_0\cap\Sigma = \{d\in\Sigma : \rho(d)=0\}$.
The remaining indices may be weakly active with zero multipliers or inactive.
If $\Sigma$ is empty, then $K=\{0\}$ and SOSC holds.
If $\mathcal A_+$ is empty, then $K_0=K$ and strong SOSC gives $Q_\tau(d)\ge d^\top W d>0$ for every nonzero $d\in K$.
Thus, for the remainder of the proof, assume that both $\Sigma$ and $\mathcal A_+$ are nonempty.
Since every term of the sum in $Q_\tau$ is nonnegative, dropping the indices outside $\mathcal A_+$ only decreases $Q_\tau$, i.e. we have

\begin{align*}
Q_\tau(d) \;\ge\; d^\top Wd + \tau\sum_{i \in\mathcal A_+}\frac{(\bar H_id)^2}{s_i^*(\tau)^2} \;=:\; \widehat Q_\tau(d).
\end{align*}

So it suffices to prove that $\widehat Q_\tau(d)>0$ on $\Sigma$ for $\tau$ small.

By the convergence assumption in the theorem, $\mu_i^*(\tau) \to \bar \mu_i>0$ for every $i\in\mathcal A_+$. The barrier complementarity equation $\mu_i^*(\tau)s_i^*(\tau)=\tau$ therefore gives, for sufficiently small $\tau>0$,

$$
s_i^*(\tau) = \frac{\tau}{\mu_i^*(\tau)} = \Theta(\tau), \qquad i \in\mathcal A_+.
$$

Let $E\coloneqq K_0\cap\Sigma$.
If $E$ is nonempty, strong SOSC and compactness imply that there is a $\sigma_E>0$ such that $d^\top W d\ge\sigma_E$ for all $d\in E$.
By continuity, choose $\delta>0$ such that $\rho(d)<\delta$ implies $d^\top W d\ge\sigma_E/2$ on $\Sigma$.
Set $\sigma\coloneqq\sigma_E/2$ in this case.
If $E$ is empty, compactness gives $\min_{d\in\Sigma}\rho(d)>0$; choose $\delta$ below this minimum and any $\sigma>0$, such that the implication $\rho(d)<\delta\implies d^\top W d\ge\sigma$ becomes trivial.

\emph{Case 1 (near $K_0$).} By this choice of $\delta$ and $\sigma$, $\widehat Q_\tau(d)\ge d^\top W d\ge\sigma>0$ whenever $d\in\Sigma$ and $\rho(d)<\delta$.

\emph{Case 2 (away from $K_0$).} On the compact set $\Sigma_\delta\coloneqq \{d\in\Sigma:\rho(d)\ge\delta\}$, boundedness of $W$ gives $M\ge0$ with $d^\top Wd\ge -M$.
We have
$$
\tau\sum_{i \in\mathcal A_+}\frac{(\bar H_id)^2}{s_i^*(\tau)^2} \ge \frac{\tau\,\rho(d)}{\max_{i \in\mathcal A_+} s_i^*(\tau)^2} \ge \frac{\tau\delta}{\max_{i \in\mathcal A_+}s_i^*(\tau)^2} = \Omega(1/\tau)\to\infty,
$$
uniformly over $\Sigma_\delta$. So there is $\bar\tau>0$ such that this sum exceeds $M+1$ for all $d\in\Sigma_\delta$, $\tau\in(0,\bar\tau]$, giving

$$
\widehat Q_\tau(d) \ge -M + (M+1) = 1 > 0, \qquad d\in\Sigma_\delta,\ \tau\in(0,\bar\tau].
$$
We have shown that $\widehat Q_\tau(d)>0$ for all $d\in\Sigma$ and $\tau\in(0,\bar\tau]$.
As $Q_\tau(d) \geq \widehat Q_\tau(d)$, the same holds for $Q_\tau$ and, by homogeneity, for all $d\in K\setminus\{0\}$. This is exactly SOSC for \eqref{eq:bqp} at $(\Delta z^*(\tau),s^*(\tau),\lambda^*(\tau),\mu^*(\tau))$ for every $\tau\in(0,\bar\tau]$.
\end{proof}

\paragraph{Example showing that original SOSC does not ensure barrier-QP SOSC.}
The following example shows that ordinary SOSC for the original NLP does not by itself ensure SOSC for the barrier QP in Theorem~\ref{thm:barrier_differentiability}.
It illustrates why the proof uses strong SOSC to establish this condition.
Consider an NLP with no equality constraints, $n_z=n_h=3$, objective
$f(z)=\tfrac12 z^\top Wz$, and inequalities $h(z)=z\le 0$.
Let
\[
e=(1,1,1)^\top,\qquad
v=\frac{(1,-1,0)^\top}{\sqrt{2}},\qquad
u=\frac{(1,1,-2)^\top}{\sqrt{6}},
\]
and define
\[
W=\frac{ee^\top}{3}-vv^\top+3uu^\top
=
\begin{bmatrix}
1/3&4/3&-2/3\\
4/3&1/3&-2/3\\
-2/3&-2/3&7/3
\end{bmatrix}.
\]
The vectors $e/\sqrt{3},v,u$ are orthonormal eigenvectors of $W$ with
eigenvalues $1,-1,3$, respectively.

\FloatBarrier
\subsection{Example of non-continuous solution map}
\label{sec:jump_ocp}
\begin{figure}[h]
	\centering
	\includegraphics[width=.49\columnwidth]{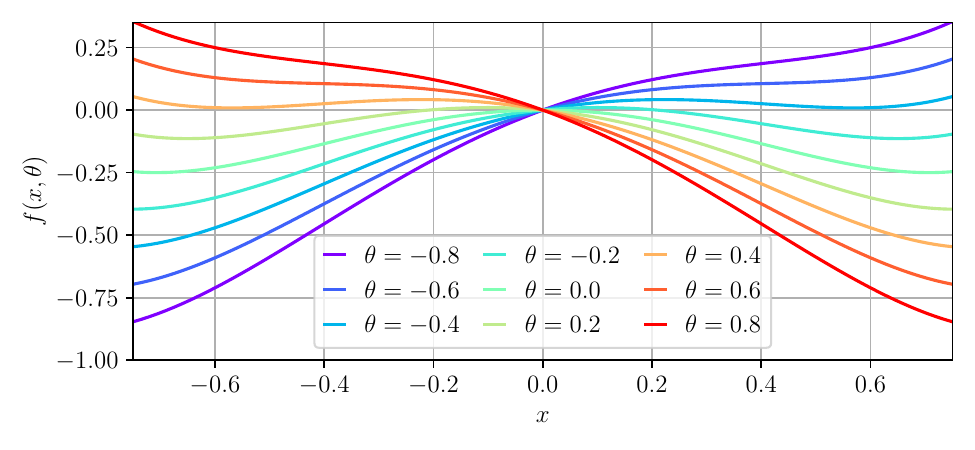}
	\includegraphics[width=.49\columnwidth]{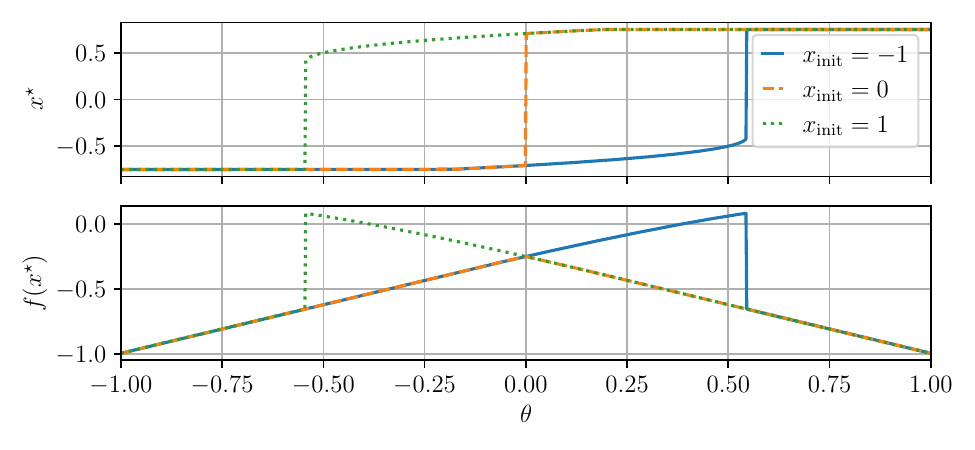}
	\caption{Example of an NLP~\eqref{eq:jump_ocp} for which the global solution $x^\star(\theta)$ jumps.
	The objective function of \eqref{eq:jump_ocp} is visualized for different values of $\theta$ on the left.
	The plots on the right show the numerical solution and corresponding objective function value obtained with different solver initializations $x\ind{init}$.
	\label{fig:jump_ocp}}
\end{figure}

Let us regard the example
\begin{mini!}
	{\substack{x}}
	{(x-1)(x+1)x^2 - \theta x}
	{\label{eq:jump_ocp}}
	{}
	\addConstraint{-0.75 \leq x}{\leq 0.75}.
\end{mini!}
For $\theta=0$, the problem has two local minimizers, with the same objective function value.
Both local minimizers exist in a neighborhood of $\theta=0$, we can denote them by $x^{\star,+}(\theta) > 0$ and $x^{\star, -}(\theta) < 0$.
For values of $\theta>0$, the local minimizer $x^{\star, +}$ is the global minimizer.
When increasing $\theta$ further, approximately at a value of 0.54 for $\theta$, the local minimizer $x^{\star, -}$ vanishes.
We visualize the solutions obtained with the NLP solver~\texttt{IPOPT}~\citep{Waechter2006} for different solver initialization $x\ind{init}$ in Figure~\ref{fig:jump_ocp}.
For $x\ind{init} = 0$, the solver converges always to the global solution, while when initializing with $x\ind{init}=-1$ or $x\ind{init}=1$ the closest local optimizer is found, which is not necessarily a global one.
At the value of $\theta$ where the local minimizer  $x^{\star, -}$ vanishes, Assumption~\ref{as:regularity}, in particular SOSC is not satisfied.

\subsection{Differentiable \acados{} layer in \texttt{leap-c}}
\label{app:leap_c}
The proposed implementation has been wrapped for integration into common ML frameworks, like \texttt{PyTorch} or \texttt{JAX}, together with additional convenience functionality to handle different kinds of parameters.
This is contained in the \texttt{leap-c} project~\citep{leap-c_releases}, which focuses on software for learning-enhanced control.
In particular, the two-solver approach and the \acados{} batch solver functionality are wrapped in the class \texttt{AcadosDiffMpcFunction}, which implements custom autograd functions in a way that is agnostic of the ML framework.
On top of that is a thin \texttt{PyTorch} specific layer called \texttt{AcadosDiffMpcTorch} which is extensively tested in RL-enhanced MPC schemes.
In addition to the solution sensitivities which are the focus of this paper, these classes also provide sensitivities of the value function, which can be leveraged for imitation learning \citep{Ghezzi2023b}.

An example of learning-enhanced MPC using this pipeline is presented in~\ref{sec:cartpole}.

\begin{figure}
	\centering
	\includegraphics[width=.45\columnwidth]{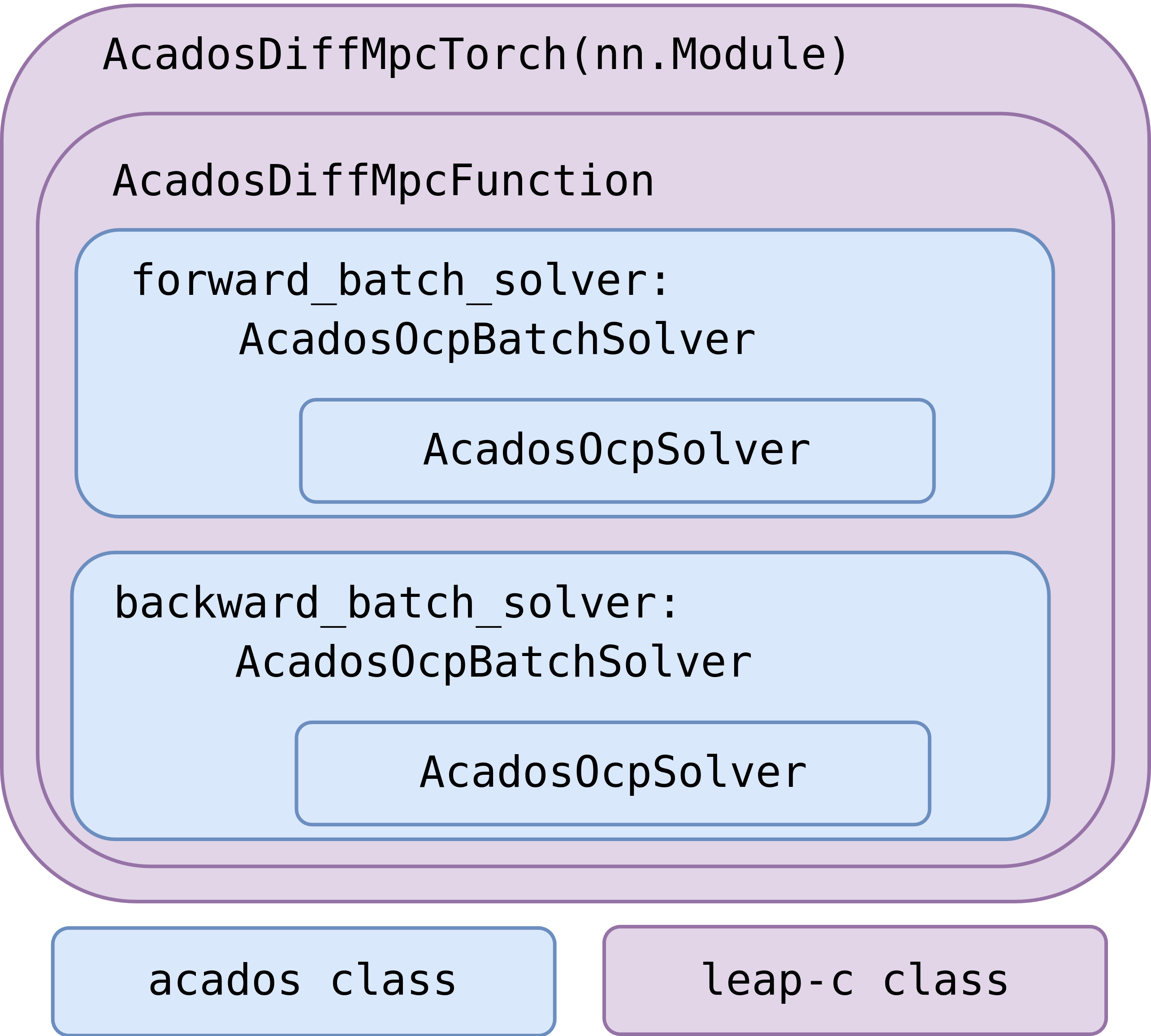}
	\caption{Schematic overview of the class hierarchy for integrating the proposed approach in common ML frameworks via \texttt{AcadosDiffMpcFunction} and specifically for \texttt{PyTorch} via \texttt{AcadosDiffMpcTorch}.}
\end{figure}

\newpage
\subsection{Usage of the proposed implementation}\label{sec:code_snippet}
The following code snippet shows how the proposed implementation can be used to compute forward or adjoint solution sensitivities.
In particular, we refer to the steps in Section~\ref{sec:sol_sens_alg} and use the two-solver approach discussed in Section~\ref{sec:acados_sol_sens}.
\begin{python}
nominal_solver.solve() # (S1)
# two-solver approach iterate transfer
iterate = nominal_solver.store_iterate()
sens_solver.load_iterate_from_obj(iterate)
# (S2), (S3)
sens_solver.setup_qp_matrices_and_factorize()
# (S4f), (S5f): forward version
result = sens_solver.eval_solution_sensitivity(stages=[0], with_respect_to='p_global')
sens_u = result['sens_u']
# (S4a), (S5a): adjoint version
s_adj = sens_solver.eval_adjoint_solution_sensitivity(seed_u=[(0, np.ones((nu, 1)))])
\end{python}

Note that the forward version computes the full Jacobian~\eqref{eq:forward_ipm_sens}, and the argument \texttt{stages = [0]} specifies that only the values corresponding to variables at stage 0 are unpacked in the \texttt{Python} wrapper.
For the adjoint version, the seed vector $\nu$ from Sec.~\ref{sec:adjoint_sens} is specified, by providing all non-zero entries, via the arguments \texttt{seed\_u} and \texttt{seed\_x}.
In the snippet above, $\nu$ is specified as a vector of zeros except for the entries corresponding to $u_0$, which are set to 1.
Since in this case all $\nu$ components corresponding to states are zero, \texttt{seed\_x} is not specified.

\subsection{Reproducibility of the presented results}\label{sec:reproducibility}
The core contribution is publicly available in the \acados{} repository~\citep{acados_releases} which is subject to the permissive 2-Clause BSD license.
The code to reproduce all Figures and Tables in this paper is available in
\anonymousOrSubmission{the zip file attached to the submission and will be publicly available in a GitHub repository and released with a zenodo doi, as done for the \acados{} releases~\citep{acados_releases}, for the final submission.}
{the public GitHub repository \url{https://github.com/FreyJo/differentiable_nmpc} with a permanant Zenodo release \url{https://doi.org/10.5281/zenodo.24015494}.}
All experiments, except for the ones corresponding to Table~\ref{tab:mcp_pytorch_gpu}, were run on a Lenovo ThinkPad T490s with an Intel Core i7-8665U CPU and 32GB of RAM running Ubuntu 22.04.

\subsection{Details on OCP with highly nonlinear parameter dependency in Example~\ref{example:solution_sens_pendulum}}
\label{app:highly_param_ocp}
We consider an OCP, which is associated with a pendulum on cart model.
The problem contains one parameter $\theta\in\R$, which enters the cost, dynamics and constraints of the OCP.
Note that the parameter does not have a physical interpretation but is introduced for illustrative purposes only.
The system is characterized by the state $x = [p, \phi, s, \omega]$, with cart position $p$, angle of the pendulum on the cart $\phi$, cart velocity $s$ and angular velocity $\omega$.
The control input $u$ is a force acting on the cart in the horizontal plane and bound to be in $[-80, 80]$.
The ODE describing the system dynamics can be found in~\citet{Verschueren2021}, with the modification that the mass of the cart $m$ is set to the parameter $\theta$.
The discrete dynamics $\phi_n$ are obtained using a fourth-order Runge-Kutta integrator with constant step size $\Delta t = \frac{T}{N}$ over a prediction horizon $T=2$ with $N=50$ shooting intervals.
The cost function is given by
\begin{align*}
	L_n(x_n, u_n; \theta) &= \theta x_n\T Q x_n + u_n\T R u_n, \\
	M(x_N; \theta) &= \theta x_N\T Q x_N
\end{align*}
with weights $Q = 2\cdot\mathrm{diag}(10^3, 10^3, 10^{-2}, 10^{-2})$ and $R = 0.2$.
In addition, we added the parametric constraint $-1.5 \leq p\theta \leq 1.5 $.
We fix the initial state to $(0, \frac{\pi}{2}, 0, 0)$ and solve the OCP for different $\theta$ values.

\subsection{Details and additional results on the benchmark problem in Example~\ref{example:mpc_pytorch_benchmark}}
\label{sec:mpc_pytorch_gpu}
\begin{table*}
\centering
\caption{Timings in [s] for solving $n_{\mathrm{batch}} \!=\! 128$ bounded LQR problems with $N \!=\! 20$, $n_x \!=\! 8$, $n_u\!=\! 4$, $n_\theta\!=\!248$.
Run on a machine with an Nvidia GeForce RTX 3080 Ti GPU and an AMD Ryzen 9 5950X 16-Core CPU utilizing the GPU capabilities of \texttt{mpc.pytorch}.
The results correspond to the earlier version of the benchmark with $A = \eye + 0.2 \cdot M$, resulting in failures for \texttt{mpc.pytorch} in case of $u\ind{max} = 1$.
\label{tab:mcp_pytorch_gpu}
}
\vspace{-1mm}
\small
\begin{tabular}{ccccccc}
\toprule
& \multicolumn{3}{c}{\textbf{Nominal solution}} & \multicolumn{3}{c}{\textbf{Solution + adjoint sens.}}\\
problem config & \acados & \texttt{mpc.pytorch} & speedup& \acados & \texttt{mpc.pytorch} & speedup\\
\midrule
$u_{\mathrm{max}} = 10^4$& 0.007 & 0.05 & 7.18 & 0.029 & 0.07 & 2.59\\
$u_{\mathrm{max}} = 1.0$& 0.008 & 15.12 & 1929.75 & 0.033 & 14.64 & 450.01\\
\bottomrule
\end{tabular}
\vspace{-2mm}
\end{table*}

The OCP with quadratic cost, linear discrete-time dynamics and bounds on the control inputs regarded in Example~\ref{example:mpc_pytorch_benchmark} can be written as
\begin{mini!}
	{\substack{x_0,\ldots, x_N, \\ u_0,\ldots, u_{N\!\shortminus 1} }}
	{\sum_{n=0}^{N\shortminus 1}
	\begin{bmatrix}
		x_n\\
		u_n
	\end{bmatrix}\T
	H
	\begin{bmatrix}
		x_n\\
		u_n
	\end{bmatrix}
	+
	x_N\T H_x x_N
	}
	{\label{eq:diff_mpc_ocp_appendix}}
	{}
	\addConstraint{x_0}{= \bar{x}_0}
	\addConstraint{x_{n+1}}{= A x_n + B u_n + b,~}{0 \leq n < N}
	\addConstraint{\!\!\!\!\!-u\ind{max}}{\leq u_n \leq u\ind{max},}{0 \leq n < N,}
\end{mini!}
where $A = \eye + 0.1 \cdot M$ and $B, b $ and $M$ consist of values sampled from a standard normal distribution.
The cost matrix $H$ is set to the identity and $H_x$ denotes the submatrix consisting of the first $\nx$ rows and columns of $H$.
The problem data~$A, B, b, H $ is regarded as parameter~$\theta$, such that $n_\theta = \nx^2+\nx\nctrl + \nx + (\nx+\nctrl)^2$.

\textbf{Benchmark design choices.}
While this paper focuses on constrained and nonlinear differentiable MPC, the benchmark problem considered is deliberately simple. This choice is primarily driven by limitations in the problem classes supported by existing methods; see Table~\ref{tab:diff_mpc}.
Restricting inequality constraints to simple control bounds facilitates a direct comparison with \texttt{mpc.pytorch}. Moreover, nonlinear problems—particularly those with nonlinear dynamics—are not supported by \texttt{cvxpygen}.

For nonlinear problems, many solvers rely on inexact or regularized Hessian approximations, which can lead to arbitrarily inaccurate sensitivities in general.
As a result, a runtime comparison would be unfair if different algorithms were allowed to employ arbitrarily coarse approximations.

Including additional software such as \texttt{diffmpc}, \texttt{trajax}, or \texttt{theseus} would require further simplification of the benchmark problem, as these tools do not support inequality constraints.
Finally, a comparison with \texttt{sIPOPT} or \texttt{CasADi}s native sensitivity approach would add limited value beyond existing studies on nominal OCP solutions, which already demonstrate that structure-exploiting OCP solvers achieve substantially lower computation times, rendering general-purpose solvers non-competitive~\citep{Vanroye2023,Verschueren2021}.

\textbf{Convergence of compared methods.}
The code accompanying this paper shows that all solvers converge to the same primal solutions and that the adjoint solution sensitivities match.

In an earlier version of this benchmark, with more difficult dynamics, using $A = \eye + 0.2 \cdot M$, and tight control bounds $u\ind{max} = 1.0$, we obtained consistent results with \texttt{acados} and \texttt{cvxpygen}, while \texttt{mpc.pytorch} failed to converge.
The convergence issues of \texttt{mpc.pytorch} can be attributed to the fact that the iLQR algorithm is based on an active-set heuristic.
In particular, it is not clear how this algorithm is supposed to remove constraints from the guess of the active set \citep{Tassa2014}.

\textbf{Evaluation metrics.}
We decided to only evaluate the pure solution and the solution together with adjoint sensitivity computation in the benchmark, as is reasoned next.
Conceptually, one can combine two methods for the nominal solution and (adjoint) sensitivity computation in a two-solver approach as outlined in Section~\ref{sec:solution_sens_nlp}.
However, when different software frameworks are used, this practically becomes nontrivial and introduces ambiguity into a corresponding runtime comparison.
Evaluating the computation time for a sole sensitivity computation, in particular for \texttt{mpc.pytorch} and \texttt{cvxpygen} is different from evaluating it together with the solution itself, as they both reuse results obtained in the nominal solve.
In particular, \texttt{mpc.pytorch} reuses the matrix factorization from the forward pass in the backward pass~\citep{Amos2018} and \texttt{cvxpygen} even uses low-rank updates to update the factorization between subsequent calls.
For these reasons, the APIs of \texttt{mpc.pytorch} and \texttt{cvxpygen} are not made for solely computing sensitivities without previously solving the problem.
In addition, the most common practical use-case is to use the same solver type for both steps, which makes regarding solution and sensitivity computation jointly meaningful for practitioners.

\textbf{Detailed discussion \& additional results.}
The speedups compared to \texttt{cvxpygen} can be attributed to different factors.
Firstly, \acados~is tailored to exploit the OCP structure, while \texttt{OSQP} and the \texttt{cvxpygen} differentiator do not exploit this special problem structure.
Secondly, the problems are solved to a tolerance of $10^{-6}$, while \texttt{OSQP} as a first order method is most suitable to achieve modest accuracies \citep{Stellato2020}.
On the other hand, this benchmark is suitable for the \texttt{cvxpygen} differentiator, since the constraints and the Hessian do not vary between the different batch instances.
This allows the \texttt{cvxpygen} differentiator to avoid a full factorization of the KKT matrix and instead perform a number of low-rank updates corresponding to the number of constraints that switch between active and inactive in subsequent batch instances, reducing the computational cost to scale quadratically in the number of variables compared to a cubic scaling for a full factorization.
In particular, there are no active inequalities for $u\ind{max} = 10^4$, and for $u\ind{max} = 1.0$ only $\approx 7.3\%$ of the inequalities are active.
For nonlinear constrained problems or nonlinear cost functions, such low-rank updates cannot be exploited, as the Hessian of the QP would always change between subsequent calls because the exact Lagrangian Hessian has to be used for sensitivity computations, as discussed in this work.

In addition to the results presented for Example~\ref{example:mpc_pytorch_benchmark}, Table~\ref{tab:mcp_pytorch_gpu} presents results obtained when utilizing the GPU capabilities of \texttt{mpc.pytorch} for the same benchmark problem on a machine with an Nvidia GeForce RTX 3080 Ti GPU and an AMD Ryzen 9 5950X 16-Core CPU.
The GPU table uses the same dimensions but the earlier dynamics $A = \eye + 0.2 \cdot M$, rather than $A = \eye + 0.1 \cdot M$ used for the current CPU benchmark.
With tight control bounds, \texttt{mpc.pytorch} fails to converge in this earlier setup.
The corresponding reported runtime ratios should therefore not be interpreted as speedups for successful solves.
The speed differences reported in Table~\ref{tab:mcp_pytorch_gpu} and Table~\ref{tab:mcp_pytorch} are specific to the hardware and problem sizes used, and the difference in dynamics also means that the CPU/GPU comparison is not controlled.
We note that the problem dimensions were chosen to represent dimensions that are typical in the context of MPC.
For the loose-bound setting, the \texttt{mpc.pytorch} timings suggest no substantial speedup from GPU use; this is consistent with the findings in \citet{Adabag2026}.

\subsection{RL-enhanced MPC for cart-pole swing-up}
\label{sec:cartpole}
\begin{figure}
	\includegraphics[width=\columnwidth]{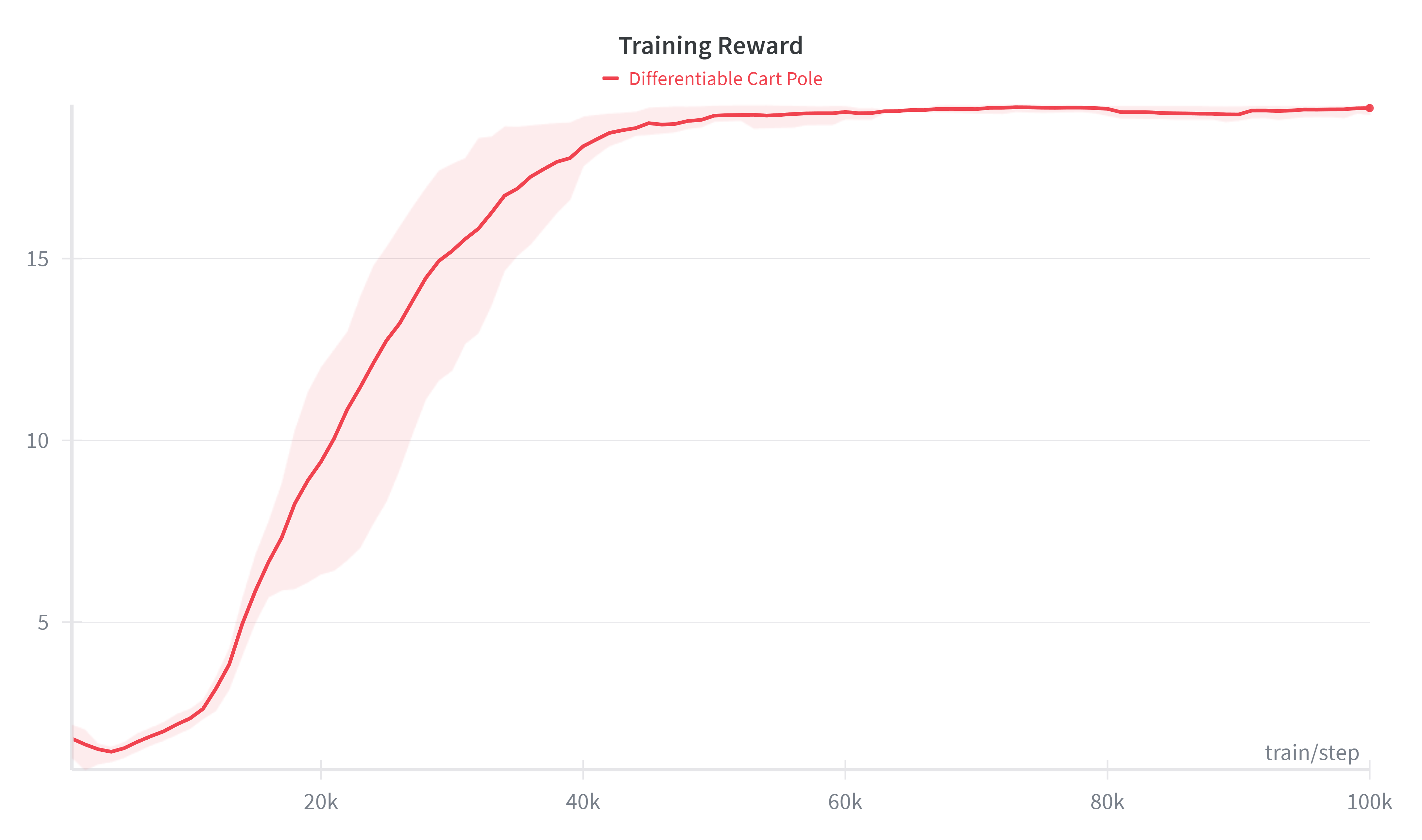}
	\caption{Training reward average and standard deviation of RL-enhanced MPC example in Sec.~\ref{sec:cartpole} over five seeds. \label{fig:leap_c_reward_up}}
\end{figure}
This section briefly presents an example of RL-enhanced MPC.
The example is implemented using \texttt{leap-c} \citep{leap-c_releases} with \texttt{acados} as the differentiable OCP solver.
The environment is a cart-pole swing-up implementation, which is based on the classic cart-pole environment from \citet{Barto1983}.
Contrary to the classic environment, the pole starts hanging down and the task is to swing it up and balance it there.

To achieve this task, the agent receives a reward that is proportional to the angular distance of the pole to the upright position.
The episode ends when the cart position exceeds the world bounds of $\pm2.4$ or after 200 steps.

The MPC controller is based on the OCP formulation
\begin{mini!}
	{\substack{x_0,\ldots, x_N, \\ u_0,\ldots, u_{N\!\shortminus 1} }}
	{\sum_{n=0}^{N\shortminus 1} L(x_n, u_n; \theta)  + M(x_N; \theta)}
	{\label{eq:cartpole_OCP}}
	{}
	\addConstraint{x_0}{= \bar x_0}
	\addConstraint{\!\!\!\!\!\!\!\!\!\!\!\!\!\!\!x_{n+1}}{=\phi(x_n, u_n; \theta), ~}{0 \leq n < N}
	\addConstraint{-80}{\leq u_n \leq 80,}{0 \leq n < N}
	\addConstraint{-2.4}{\leq x_{1,n} \leq 2.4, \label{eq:leap_c_xbound}}{1 \leq n < N.}
\end{mini!}
The dynamics function $\phi$ represents the discretized ODE of the system for a time interval of $0.25\mathrm{s}$, which coincides with one step in the environment.
The stage and terminal cost functions $L$, $M$ are of linear least-squares type with fixed weighting matrices, a learned reference value for the angle, and reference values zero for all other states and the control.
The box constraints correspond to actuator bounds and the world bounds on the cart position.
The horizon of the MPC controller is set to $N=5$.

To learn the reference parameter, a variant of Soft Actor-Critic \citep{Haarnoja2018} from \texttt{leap-c} is used, which in \texttt{leap-c} is called \texttt{sac-fop}.
Instead of a Multilayer Perceptron (MLP) generating actions directly like in SAC, the actor in this variant consists of an MLP generating parameters for the MPC controller, which then computes the action.
While training, the backpropagation uses the solution sensitivities to differentiate through the MPC controller, essentially treating the MPC controller as a final layer inside the actor, similar to~\citet{Romero2024}.

Figure~\ref{fig:leap_c_reward_up} shows the evolution of the reward through the training.
The reward is scaled such that one step can achieve a maximum of 0.1 reward, resulting in a maximum cumulative reward of 20 over one episode, hence it can be seen that the agent learns to successfully balance the pole.
Note that due to the short horizon, even with the reference angle fixed to 0, the nominal MPC controller cannot balance the pole properly.
This demonstrates that the proposed approach can be effectively integrated into a learning pipeline and achieves substantially better performance than MPC alone.

The results can be reproduced by running the script in the \texttt{leap-c} repository\footnote{\url{https://github.com/leap-c/mpc-sac/blob/d38726d76c9477274f7a588fe27a24e40120d897/scripts/run_sac_fop.py}}.

\section*{Checklist}

\begin{enumerate}

\item {\bf Claims}
    \item[] Question: Do the main claims made in the abstract and introduction accurately reflect the paper's contributions and scope?
    \item[] Answer: \answerYes{} %
    \item[] Justification: The abstract and Introduction state three contributions, each supported by the analysis, implementation, and experiments cited below.
    \begin{itemize}
		\item \textbf{Claim 1:} ``First, this work discusses the efficient computation of parametric (adjoint) solution sensitivities for general nonlinear programs (NLP) via the implicit function theorem (IFT).
		In particular, this work highlights common pitfalls of the IFT and discusses their mitigation.
		Previous work in the machine learning (ML) literature was subject to this pitfall, or limited to convex or unconstrained nonlinear least squares problems.''\\
        \textbf{Justification:} Section~\ref{sec:preliminaries} states the regularity assumptions and develops the IFT-based sensitivity analysis for general NLPs.
        Remark~\ref{remark:hess_approx} and Example~\ref{example:solution_sens_pendulum} show how an approximate Hessian can degrade sensitivities.
        The related-work discussion and Table~\ref{tab:diff_mpc} support the stated comparison with prior methods.

		\item \textbf{Claim 2:} ``Second, we introduce SBQP, which generalizes SQP with an IPM QP solver by allowing the barrier parameter to approach a prescribed nonnegative target value.
		We detail how IPMs or SBQP can be used to obtain a differentiable approximation of the solution map and its sensitivities suitable for gradient-based learning algorithms by leveraging the smoothing of the optimality conditions that these methods employ.''\\
        \textbf{Justification:} Section~\ref{sec:sqp-ipm} defines SBQP and its zero-target case.
        Theorems~\ref{thm:nlp-sens-exist}, \ref{thm:smooth_KKT_sol}, and~\ref{thm:barrier-bqp} in Section~\ref{sec:dense_solution_sens} establish the local sensitivity and smoothing results under their stated assumptions.
        Section~\ref{sec:dense_solution_sens} discusses the relevance of smoothed sensitivities to learning, and Section~\ref{sec:cartpole} demonstrates their use in an RL example.

		\item \textbf{Claim 3:} ``Finally, we present an efficient implementation of a differentiable solver for optimal control structured NLPs within the open-source software package \acados{} which provides both forward and adjoint sensitivities.
		In particular, the optimal control problem (OCP) formulation covers general, potentially nonlinear, costs and constraints on both states and controls while allowing all problem functions, i.e., costs, constraints and dynamics, to be fully parametric.
		A comparison with state-of-the-art solvers demonstrates speedup factors greater than three on a CPU.''\\
        \textbf{Justification:} Sections~\ref{sec:diff_mpc} and~\ref{sec:acados_sol_sens} describe the implementation, and Examples~\ref{example:solution_sens_pendulum} and~\ref{example:chain_sens} demonstrate parametric nonlinear OCPs and sensitivity computation.
        Example~\ref{example:mpc_pytorch_benchmark} and Table~\ref{tab:mcp_pytorch} support the CPU speedup claim for the evaluated linear-quadratic OCP benchmark with simple input bounds; the reported speedups are specific to this benchmark.
    \end{itemize}

\item {\bf Limitations}
    \item[] Question: Does the paper discuss the limitations of the work performed by the authors?
    \item[] Answer: \answerYes{} %
    \item[] Justification: The limitations of the approach presented in the paper can be structured into fundamental theoretical ones and ones that concern the proposed implementation.\\
    The theoretical limitations that need to hold for the existence of the solution sensitivities of NLP~\ref{eq:z_star_def} are given in Assumption~\ref{as:regularity}.
	Example~\ref{example:solution_sens_tutorial} showcases that this assumption is not satisfied at active-set changes, but that reasonable results can be obtained by the smoothing employed by an interior-point method.
	The example in Appendix~\ref{sec:jump_ocp} shows the fundamental issue of general nonconvex NLPs with multiple local minima and solution maps that exhibit jumps.
	Their treatment in practice is not part of this work and both examples violate the conditions in Assumption~\ref{as:regularity}.\\
	The limitations of the proposed implementation are given in the paragraph ``Limitations'' in Section~\ref{sec:acados_sol_sens}.

\item {\bf Theory assumptions and proofs}
    \item[] Question: For each theoretical result, does the paper provide the full set of assumptions and a complete (and correct) proof?
    \item[] Answer: \answerYes{} %
	\item[] Justification: The main text presents three theorems in Section~\ref{sec:dense_solution_sens}; the Appendix also states the IFT and the NLP/QP sensitivity result.
	The assumptions are stated alongside the results or referenced through Assumption~\ref{as:regularity}.
	The Appendix proves Theorem~\ref{thm:nlp-sens-exist}, Theorem~\ref{thm:barrier-bqp}, and Theorem~\ref{thm:nlp-qp}.
	Theorem~\ref{thm:smooth_KKT_sol} follows from Properties 1 and 3 in~\citet{Pirnay2012}, and a comprehensive proof is given alongside Theorem 8 in~\citet{Wright1992}.
	The IFT is stated in the Appendix in its classical form and cited to~\citet{Dini1907}.

    \item {\bf Experimental result reproducibility}
    \item[] Question: Does the paper fully disclose all the information needed to reproduce the main experimental results of the paper to the extent that it affects the main claims and/or conclusions of the paper (regardless of whether the code and data are provided or not)?
    \item[] Answer: \answerYes{} %
    \item[] Justification: All information needed to reproduce the figures and tables in this paper is contained in the code in the supplemental material as described in Section~\ref{sec:reproducibility}.
    \item[] Guidelines:

\item {\bf Open access to data and code}
    \item[] Question: Does the paper provide open access to the data and code, with sufficient instructions to faithfully reproduce the main experimental results, as described in supplemental material?
    \item[] Answer: \answerYes{} %
    \item[] Justification: All code to reproduce the figures and tables in this paper is available in the supplemental material as described in Section~\ref{sec:reproducibility}.

\item {\bf Experimental setting/details}
    \item[] Question: Does the paper specify all the training and test details (e.g., data splits, hyperparameters, how they were chosen, type of optimizer, etc.) necessary to understand the results?
    \item[] Answer: \answerYes{} %
    \item[] Justification: All code, (including test details, parameter choices etc.) to reproduce the figures and tables in this paper is available in the supplemental material as described in Section~\ref{sec:reproducibility}.

\item {\bf Experiment statistical significance}
    \item[] Question: Does the paper report error bars suitably and correctly defined or other appropriate information about the statistical significance of the experiments?
    \item[] Answer: \answerNA{} %
    \item[] Justification: The paper does not contain experiments for which this question is relevant.

\item {\bf Experiments compute resources}
    \item[] Question: For each experiment, does the paper provide sufficient information on the computer resources (type of compute workers, memory, time of execution) needed to reproduce the experiments?
    \item[] Answer: \answerYes{} %
    \item[] Justification: The experiments presented in this paper can be run on a modern laptop.
	The runtime of the experiment in Example~\ref{example:mpc_pytorch_benchmark} is roughly 30 minutes.
    The runtime for the Examples~\ref{example:chain_sens} and~\ref{example:solution_sens_pendulum} is below 5 minutes.

	The runtime for the example in Section~\ref{sec:cartpole} is around 10 hours on a modern laptop.
	This example only shows how the building block of differentiable NMPC can be used in a learning framework, while this is not the focus of this work.
    
\item {\bf Code of ethics}
    \item[] Question: Does the research conducted in the paper conform, in every respect, with the NeurIPS Code of Ethics \url{https://neurips.cc/public/EthicsGuidelines}?
    \item[] Answer: \answerYes{} %
    \item[] Justification: After carefully reading the guidelines, the authors do not see any conflicts with the NeurIPS Code of Ethics.

\item {\bf Broader impacts}
    \item[] Question: Does the paper discuss both potential positive societal impacts and negative societal impacts of the work performed?
    \item[] Answer: \answerYes{} %
    \item[] Justification:
	This work has the potential of advancing learning-enhanced MPC by combining explicit constraint handling with solution sensitivity computation.
	As previous related works mostly relied on heuristic constraint treatment, this can reduce operational risks and ensure that controllers respect physical and safety-critical limits.
	While the methods may be applied broadly, including in contexts without direct societal benefit, their impact ultimately depends on responsible deployment, highlighting the dual-use nature of the technology.

\item {\bf Safeguards}
    \item[] Question: Does the paper describe safeguards that have been put in place for responsible release of data or models that have a high risk for misuse (e.g., pretrained language models, image generators, or scraped datasets)?
    \item[] Answer: \answerNA{} %
    \item[] Justification: The answer NA means that the paper poses no such risks.

\item {\bf Licenses for existing assets}
    \item[] Question: Are the creators or original owners of assets (e.g., code, data, models), used in the paper, properly credited and are the license and terms of use explicitly mentioned and properly respected?
    \item[] Answer: \answerYes{} %
    \item[] Justification: The paper only uses open-source software as dependencies. All assets used are properly credited to the best of our knowledge.

\item {\bf New assets}
    \item[] Question: Are new assets introduced in the paper well documented and is the documentation provided alongside the assets?
    \item[] Answer: \answerYes{} %
    \item[] Justification: The core software contribution is publicly available in the~\acados~repository~\cite{acados_releases} which is subject to the permissive 2-Clause BSD license.
	Its documentation is publicly available on \url{https://docs.acados.org/}.
	In particular, the functions introduced and showcased in this work, see Section~\ref{sec:code_snippet} are contained in the documentation.

\item {\bf Crowdsourcing and research with human subjects}
    \item[] Question: For crowdsourcing experiments and research with human subjects, does the paper include the full text of instructions given to participants and screenshots, if applicable, as well as details about compensation (if any)? 
    \item[] Answer: \answerNA{} %
    \item[] Justification: The answer NA means that the paper does not involve crowdsourcing nor research with human subjects.

\item {\bf Institutional review board (IRB) approvals or equivalent for research with human subjects}
    \item[] Question: Does the paper describe potential risks incurred by study participants, whether such risks were disclosed to the subjects, and whether Institutional Review Board (IRB) approvals (or an equivalent approval/review based on the requirements of your country or institution) were obtained?
    \item[] Answer: \answerNA{} %
    \item[] Justification: The answer NA means that the paper does not involve crowdsourcing nor research with human subjects

\item {\bf Declaration of LLM usage}
    \item[] Question: Does the paper describe the usage of LLMs if it is an important, original, or non-standard component of the core methods in this research? Note that if the LLM is used only for writing, editing, or formatting purposes and does \emph{not} impact the core methodology, scientific rigor, or originality of the research, declaration is not required.
    \item[] 
    \item[] Answer: \answerNA{} %
    \item[] Justification: The answer NA means that the core method development in this research does not involve LLMs as any important, original, or non-standard components.

\end{enumerate}

\end{document}